\documentclass[12pt]{article}
\usepackage{graphicx}
\usepackage{graphics,amsmath,amssymb}

\catcode`@=11
\@addtoreset{equation}{section}
\catcode`@=12

\newtheorem{theo}{Theorem}
\newtheorem{rem}{Remark}[section]
\newtheorem{lem}{Lemma}[section]
\newtheorem{conj}{Conjecture}
\newtheorem{prop}{Proposition}[section]
\newtheorem{cor}{Corollary}[section]
\newtheorem{dfn}{Definition}[section]

\newcommand{\mN}{\mathbb{N}}

\newcommand{\mR}{\mathbb{R}}
\newcommand{\mT}{\mathbb{T}}
\newcommand{\mZ}{\mathbb{Z}}

\newcommand\thet\vartheta
\newcommand\eps\varepsilon
\newcommand\ph\varphi
\newcommand\kap\varkappa

\newcommand\bG {\mbox{\bf G}}

\newcommand\bH {\mbox{\bf H}}

\newcommand\bO {\mbox{\bf O}}
\newcommand\br {\mbox{\bf r}}
\newcommand{\calA}{\mathcal{A}}

\newcommand\bU {\mbox{\bf U}}

\newcommand\SM {{\cal S\!M}} 
\newcommand\BB {{\cal B}}

\newcommand\DD {{\cal D}}
\newcommand\EE {{\cal E}}

\newcommand\JJ {{\cal J}}
\newcommand\KK {{\cal K}}

\newcommand\NN {{\cal N}}
\newcommand\OO {{\cal O}}
\newcommand\PP {{\cal P}}

\newcommand\TT {{\cal T}}

\newcommand\WW {{\cal W}}

\newcommand\diam {\operatorname{diam}}
\newcommand\dist {\operatorname{dist}}
\newcommand\grad {\operatorname{grad}}

\newcommand\Int {\operatorname{int}}

\newcommand\spn {\operatorname{span}}

\newcommand\vol {\operatorname {vol}}
\newcommand{\eqdef}{\stackrel{\mathrm{def}}{=}}

\newcommand\qed{{\unskip\nobreak\hfil\penalty50
 \hskip2em\hbox{}\nobreak\hfil\mbox{\rule{1ex}{1ex} \qquad}
   \parfillskip=0pt \finalhyphendemerits=0\par\medskip}}

\begin{document}
\title{Arnold diffusion in multidimensional a priori unstable Hamiltonian systems}

\author{
M. Davletshin\\
Independent University of Moscow\\
D.~Treschev\\
Steklov Mathematical Institute and Moscow State University
}
\maketitle

\begin{abstract}
We study the Arnold diffusion in a priori unstable near-integrable systems in a neighbourhood of a resonance of low order. We consider a non-autonomous near-integrable Hamiltonian system with $n+1/2$ degrees of freedom, $n\ge 2$. Let the Hamilton function $H$ of depend on the parameter $\eps$, for $\eps=0$ the system is integrable and has a homoclinic asymptotic manifold $\Gamma$. Our main result is that for small generic perturbation in an $\eps$-neighborhood of $\Gamma$ there exist trajectories the projections of which on the space of actions cross the resonance. By ``generic perturbations'' we mean an open dense set in the space of $C^r$-smooth functions $\frac{d}{d\eps}\big|_{\eps=0} H$, $r=r_0,r_0+1,\ldots,\infty,\omega$. Combination of this result with results of \cite{DT} answers the main questions on the Arnold diffusion in a priori unstable case: the diffusion takes place for generic perturbation, diffusion trajectories can go along any smooth curve in the action space with average velocity of order $\eps/|\log \eps|$.
\end{abstract}

\section{Introduction}
Arnold in \cite{Arn64} proposed an example of a near-integrable
Hamiltonian system
\begin{eqnarray}
\label{as}
& \dot x = \partial H / \partial y, \quad
  \dot y = - \partial H / \partial x, \quad
  H = H_0(y) + \eps H_1(x,y,t) + O(\eps^2),  & \\
\nonumber
& x\in\mT^n = \mR^n/\mZ^n, \quad
  y\in\mR^n, \quad
  t\in\mT,  \quad
  0\le\eps\ll 1 &
\end{eqnarray}
with the convex in the actions $y$ Hamiltonian $H_0$ and $n=2$, where the slow variables $y$ can change by a quantity of order 1 on a trajectory for all sufficiently small $\eps>0$. Numerical simulations show that evolution of $y$ is a combination of small oscillations with more or less random drift.
Therefore, Chirikov proposed to call this phenomenon the Arnold diffusion.
A clear explanation and the geometrical mechanism of the diffusion is given in \cite{KL}.

The main problem associated to the Arnold diffusion in systems (\ref{as}) is its genericity and in what sense this genericity is understood.
First of all this depends on the functional space which $H$ belongs to.
Physically the most interesting case is real-analytic.
According to the Nekhoroshev theory \cite{Nekh}, in real-analytic systems, satisfying the so-called steepness condition, average velocity of the action drift along a trajectory is estimated from above by an exponentially small quantity $\exp(-\alpha\eps^{-\beta})$ with positive $\alpha$ and $\beta$.
For more details, see \cite{AKN}.
Such exponentially small effects make the problem of genericity in real-analytic case very difficult and now there is no clear ways to its solution. The smooth case is much simpler to analyze but it also contains complicated technical problems.
It is natural to expect that the following conjecture is true:

\begin{conj}
\label{conj1}

{\bf A.} The diffusion exists for a typical set of perturbations.

{\bf B.} Evolution of the slow variables occurs along any smooth curve in the space of slow variables.
\end{conj}

Essential progress in formalisation and proof of it appeared in recent preprints \cite{GuaKalZh} and \cite{KZZ}.
Kaloshin and Zhang in \cite{Kal-Zhang5,Kal-Zhang2} (see also Marco \cite{Marco})
prove genericity of the diffusion in systems with $2.5$ degrees of freedom. In \cite{Kal-Zhang3} a proof in systems with $3.5$ degrees of freedom is announced. In \cite{Kal-Zhang4} an autonomous case $n\ge 4$ is considered and methods that can be used in this situation are discussed. In \cite{Cheng} a proof for $n \ge 4$ is announced (see also \cite{ChengMRL, ChengCJM, ChengJDG}).

There are several simpler situations, where Arnold diffusion occurs without exponentially small effects.
One of them is the diffusion in the so-called a priori unstable systems.
The Hamiltonian of a near-integrable non-autonomous a priori unstable system is as follows:
\begin{equation}
\label{eq:ham_ini}
    H(y,x,v,u,t,\eps)
  = H_0(y,v,u) + \eps H_1(y,x,v,u,t)+O(\eps^2).
\end{equation}
The dynamics of the system with unperturbed Hamiltonian $H_0$ is a combination of a system with one degree of freedom with hyperbolic equilibrium state and an $n$-dimensional rotator. Here $(v,u)\in D\subset \mR^2$, $y=(y_1,\ldots,y_n)\in\overline\DD\subset\mR^n$, $x=(x_1,\dots,x_n)\in\mT^n$, $\eps\ge 0$ is a small parameter, where $\DD \subset \mR^n$ is an open domain with compact closure $\overline\DD$. The symplectic structure is $dx\wedge dy+du\wedge dv$. The Hamiltonian (\ref{eq:ham_ini}) is time-periodic with period $1$ and satisfies several conditions:

${\bf H_0 1}$. {\it  In the unperturbed Hamiltonian the variables $y$ are separated from $u$ and $v$ i.e., $H_0(y,v,u) = F(y,f(v,u))$.}

${\bf H_0 2}$. {\it The function $f$ has a nondegenerate saddle point $(v,u) = (0,0)$. This is a unique critical point on a compact connected component of the set
$$
  \gamma = \{(v,u)\in D : f(v,u) = f(0,0)\} .
$$}
In other words, $(0,0)$ is a hyperbolic equilibrium state of the Hamiltonian system $(D,dv\wedge du,f)$ with one degree of freedom.

We put
$$
  E(y) = H_0(y,0,0),\quad
  \hat \nu = \partial E/\partial y : \overline\DD \to \mR^n.
$$

${\bf H_0 3}$. {\it For any $y\in\overline\DD$ the
$(n\times n)$-matrix $\partial^2 E / \partial y^2$ is nondegenerate. This means that the map $y \mapsto \hat \nu(y)$ is locally invertible in a neighbourhood of any $y_0 \in \DD$. }
\medskip

${\bf H_1 1}$. {\it  $H_1(y,x,v,u,t) \in C^{\bf r}(\DD\times\mT^n\times D\times\mT)$, where ${\bf r}>n+3$.}\\

The study of the diffusion (i.e. a drift of the slow variables $y$) in a priori unstable systems deals with three aspects presened in the following conjecture \cite{AKN}:
\begin{conj}
\label{conj:diffu}

{\bf A.} The diffusion exists for an open and dense set of $C^{\mathbf r}$-perturbations. % $H_1$.

{\bf B.} Evolution of the slow variables $y$ occurs along any smooth curve $\chi\subset\DD$.

{\bf C.} There are ``fast'' diffusion trajectories whose average velocity along $\chi$ is of order $\eps/|\log\eps|$.
\end{conj}
For $n=1$ Conjecture \ref{conj:diffu} was proved in  \cite{Tre_au}. For $n>1$ there are only partial results.

There are several approaches to solution of these problems. The first one is the construction of hyperbolic tori which form  transition chains with heteroclinic connections \cite{Arn64,CG,DLS2,Fontich,GGM}. Later this approach was supplemented with ideas of the scattering map
\cite{DLS_scat,DLS,Delsh-Hug1,Delsh-Hug2} and symbolic dynamics \cite{Bounemoura}.
In \cite{DLS} Arnold diffusion in the case of several hyperbolic directions (in our notations, the case of vector variables $u,v$) is studied.

Another approach uses variational methods: \cite{Bessi,Berti,BeBo}. In \cite{Cheng-diff1} the existence of Arnold diffusion is proved for systems with 2.5 degrees of freedom for an open and dense set of perturbations and convex $H_0$.
In \cite{Cheng-diff2} it is proved for multidimensional ``hyperbolic'' variables $u,v$.

If the unperturbed system is non-integrable (a priori chaotic case) then diffusion gets faster and more pronounced, however its general mechanisms remain essentially the same:
\cite{BT99,DLS00,Pif,DLS06,GT08,GT17}

Also note two following papers.
In the preprint \cite{GuaKalZh} the formulas for the separatrix map are improved up to terms of the second order when the perturbation is a trigonometric polynomial in angle variables and time. Using this result, in \cite{KZZ} the existence of so-called normal hyperbolic invariant foliations for an open set of trigonometric polynomial is proven.
This result is used in \cite{KZZ} to prove the existence of the so-called normally hyperbolic invariant foliations for an open set of trigonometric perturbations. For $\eps \to 0$ after a proper scaling, restriction of the dynamics to these invariant foliations turns out to be equivalent to certain diffusion process.

In this paper we continue research of Arnold diffusion in a priori unstable systems that was started in \cite{Tre_au,DT,DavTre16}.
We use methods of the multidimensional separatrix map that have been developed in \cite{Tre_sm_gl,Tre_sm_sym,PifTre} and the idea of the anti-integrable limit \cite{aubry0, BT2}.

According to ${\bf H_0 2}$ the separatrices $\gamma$ are doubled in the one degree of freedom system with Hamiltonian $f$. These separatrices are homeomorphic to a figure-eight: two loops, $\gamma^\pm$, issuing from one point,
$\gamma = \gamma^+ \cup \gamma^-$. The phase flow of the system generates orientation on $\gamma^\pm$.
The orientation on $D$ is determined by the orientation of the coordinate system $u,v$.
Without loss of generality we will assume that
the orientation of $\gamma^\pm$ coincides with the orientation of the domain $D$, i.e. the motion on the separatrices is counterclockwise.

The frequency vector $ \nu=(-\hat \nu,1)$ is called resonant if there exists a non-zero
vector $k =(\hat k,k_0)\in \mZ^{n+1}$ such that $\langle  \nu,k\rangle=0$. Below we will call the vector $k$ also resonant.
This resonance determines the hypersurface
\begin{equation}
S_0^k=\{y\in \overline \DD \colon \langle  \nu(y),k\rangle=0 \}.
\end{equation}
If $|k| \leq
C$, where $C$ does not depend on $\eps$, then this resonance is called $C$-strong or strong.
Strong resonances split the domain $\overline \DD$ into finite number of connected components.
The points of intersection of two resonant hypersurfaces that correspond to non-collinear integer vectors $k_1$ and $k_2$ are called multiple resonances. If $|k_i|\leq C$, $i=1,2$, then we call such points multiple strong resonance.

Now we describe the main result of \cite{DT}. For $\delta\ge 0$ and $k \in\mZ^{n+1}\setminus\{0\}$ define a neghbourhood of the resonance in $\overline \DD$
\begin{equation}
\label{StildeS}
    S^k_{\delta}
  = \{\eta\in\overline\DD : |\langle k, \nu(\eta)\rangle|
               \le \delta\}.
\end{equation}
Consider the set  $Q=\overline \DD \setminus \cup_{0 < |k| \le C_{\diamond}} S_{\delta}^k$ for
\begin{equation}
\label{delta}
\delta=O\big(|\log^{-1} \eps|\big).
\end{equation}

\begin{theo} \cite{DT}
\label{DTtheo}
For an open dense set in the $C^{\bf r}$ space of functions $H_1$ there exists a constant
$C_{\diamond}=C_{\diamond}(H_0,H_1)$ that does not depend on $\eps$ and there exist constants
$\alpha=\alpha(\br),\eps_0,c_d,c_v > 0$ such that for every smooth curve $\chi\subset Q$ with
endpoints $\chi_0,\chi_1$ and length $|\chi|$, and for all positive $\eps<\eps_0$, the perturbed
system has a trajectory
\begin{equation}
\label{dif_traj}
  \big( y(t), x(t), v(t), u(t) \big), \quad
  t\in [0,T]
\end{equation}
with the following properties:

(i) $|y(0) - \chi_0| < c_d\,|\log\eps|^\alpha \eps^{1/4}$,
    $|y(T) - \chi_1| < c_d\,|\log\eps|^\alpha \eps^{1/4}$,

(ii) the curve $\{y(t) : t\in [0,T]\}$ lies in the
$c_d\,|\log\eps|^\alpha \eps^{1/4}$-neighbourhood of $\chi$,

(iii) $c_v T \eps / |\log\eps| < |\chi|$.
\end{theo}
Thus, in \cite{DT} Conjecture \ref{conj:diffu} is proved {\it far} from strong resonances.
In \cite{DavTre16} the existence of the diffusion in the neighbourhood of strong resonances was proved
for $H_1$ that are generic trigonometric polynomials in angles $x$.

In this paper we prove Conjecture \ref{conj:diffu} for $n> 1$ in the vicinity of strong resonances.
We construct a diffusion trajectory that crosses the neighbourhood $S^k_{\delta}$ of a $C_{\diamond}$-strong resonance\footnote{below for simplicity we call it the strong resonance}, where the constant $C_{\diamond}$ is from Theorem \ref{DTtheo}, and obtain an estimate for the time of such transition.
The main difficulty of the diffusion problem near resonances of a low order consists in the fact that, generally speaking, we cannot construct a trajectory of the separatrix map that monotonously approaches the resonance.

Thus, we consider non-autonomous Hamiltonian systems with Hamiltonian (\ref{eq:ham_ini}) satisfying the conditions $\bf H_01-\bf H_03$.
Now we formulate the main result of the paper.

Let $\chi \subset \DD$ be a piecewisely smooth curve with the endpoints $y_0$, $y_1$ which lie on the opposite components of the boundary of $S^k_{\delta}$.
Without loss of generality
\begin{equation}
\label{bound}
\langle \hat k,\hat\nu(y_0)\rangle-k_0=\delta, \quad \langle \hat k,\hat\nu(y_1)\rangle-k_0=-\delta.
\end{equation}
Moreover, let the curve $\chi$ and $S_0^k$ transversally intersect at a point $y_*$ which is not a multiple resonance.
Then a small neighbourhood of this point does not contain any other strong resonance $S_0^l$ with $l\nparallel k$. This neighbourhood does not depend on $\eps$.
In our situation the set of multiple strong resonances is a finite union of manifolds of codimensions not less than 2.

Denote by
$$
O(\chi,r)\eqdef\bigcup_{y \in \chi}B(y,r).
$$
the neighbourhood of the curve $\chi$.
Here $B(y,r)$ is a ball of a center $y$ and radius $r$.
\begin{theo}
\label{main}
Let the function $H_0$ satisfy conditions $\bf H_01$--$\bf H_03$.
Then for any generic $H_1$ and some constants $\eps_0(H_0,H_1), c_v>0$ for all $\eps \in (0,\eps_0)$ the perturbed system has a trajectory
 \begin{equation}
 \label{traj}
 (y(t),x(t),v(t),u(t)), \qquad t\in[0,T]
 \end{equation}
such that:

(i) $|y(0)-y_0|<\eps^{1/(7n)}$, $|y(T)-y_1|<\eps^{1/(7n)}$,

(ii) the projection of (\ref{traj}) to $\DD$ lies in $O(\chi,\eps^{1/(8n)})$,

(iii) $c_vT\eps < 1$.
\end{theo}

Now we explain what we mean by generic or typical perturbation.
\begin{dfn}
\label{typ:dfn}
A subset $X \subset C^{\bf r}(\overline \DD\times\mT^n\times D\times \mT)$ is called typical, if two conditions hold:

1) $X$ is open and dense in $C^{\bf r}$-topology,

2) $X$ forms a set of full measure in any typical in the first sense family of functions in $C^{\bf r}(\overline \DD\times\mT^n\times D\times \mT)$ smoothly depending on one or several parameters.
\end{dfn}
Below we formulate explicit conditions for the perturbations $H_1$: see  $\bf H_1 1$--$\bf H_1 4$.

Combining Theorem \ref{DTtheo} in the domain $Q$ and Theorem \ref{main} in a neighbourhood of strong resonances, we obtain a proof of Conjecture \ref{conj:diffu} in the whole domain $\DD$. More precisely, as a corollary we have
\begin{theo}
\label{everywhere}
Suppose that the functions $H_0$ and $H_1$ satisfy the conditions $\bf H_01$--$\bf H_03$, $\bf H_11$--$\bf H_14$.
Let $\chi \subset \DD$ be a smooth curve of a finite length $|\chi|$ with endpoints $\chi_0$, $\chi_1$ that transversally crosses resonance hypersurfaces $S^k_0$, $0 < |k| \le C_{\diamond}(H_0,H_1)$ and does not contain multiple strong resonances.
Then there exists a constant  $c_v>0$ such that for sufficiently small $\eps>0$ the perturbed system has a trajectory (\ref{traj})
with the following properties:

(i) $|y(0)-\chi_0|<\eps^{1/(7n)}$, $|y(T)-\chi_1|<\eps^{1/(7n)}$,

(ii) the projection of (\ref{traj}) onto $\DD$ lies in $O(\chi,\eps^{1/(8n)})$,

(iii) $c_vT\eps/|\log\eps| < |\chi|$.
\end{theo}

Theorem \ref{everywhere} implies a similar theorem for an arbitrary piecewisely smooth curve $\chi \subset \DD$ of a finite length. But in condition $(i)$ instead of $\eps^{1/(7n)}$ we will have to use an arbitrary small constant $\epsilon$ which is independent of $\eps$. Since multiple strong resonances form the set of codimension not less than $2$, we can replace the curve $\chi$ by a curve $\tilde \chi$ that is $\epsilon$-close to $\chi$ and satisfies conditions of Theorem  \ref{everywhere}.
%\begin{theo}
%\label{everywhere2}
%Let $\chi \subset \DD$ be an arbitrary piecewisely smooth curve of a finite length with endpoints $\chi_0$, $\chi_1$.
%Suppose also that the functions
%$H_0$ and $H_1$ satisfy the conditions $\bf H_01-\bf H_03$, $\bf H_11-H_14$. Then for every $\epsilon >0$ there exist constants  $c_v>0$, $\eps_0>0$ such that for all $\eps \in (0,\eps_0)$ the perturbed system has a trajectory
% \begin{equation}
% (y(t),x(t),v(t),u(t)), \qquad t\in[0,T],
% \end{equation}
%with the following properties:

%(i) $|y(0)-\chi_0|<\epsilon$, $|y(T)-\chi_1|<\epsilon$,

%(ii) the projection (\ref{traj}) onto $\DD$ lies in $O(\chi,\epsilon)$,

%(iii) $c_vT\eps/|\log\eps| < |\chi|$.
%\end{theo}
\subsection{Structure of the paper}
The unperturbed system has the normally hyperbolic manifold $N = \{u=v=0\}$, foliated by invariant tori
\begin{equation}
\label{N(y)}
  N(y^0) = \{y = y^0, u=v=0\}.
\end{equation}
Consider the asymptotic manifolds
\begin{equation}
\label{Gamma}
   \Gamma^\pm
 = \bigcup_{y\in\overline\DD}
       \{y\}\times\mT^n\times\gamma^\pm\times\mT.
\end{equation}
The diffusion trajectory will be constructed in a small neighborhood of the set
$\Gamma = \Gamma^+ \cup \Gamma^-$.
More precisely, after the perturbation $N$ deforms a little and remains a nomally hyperbolic manifold $N_{\eps}$. The corresponding asymptotic manifolds $\Gamma^{\pm}_{\eps}$ in general no more coincide and split. The diffusion trajectories we construct are contained in a $c\eps$-neighbourhood of $\Gamma^+_{\eps} \cup \Gamma^-_{\eps}$ where $c$ is a small constant. In fact, we do not use the fact of the existence of $N_{\eps}$ and $\Gamma^{\pm}_{\eps}$. Our method is based on the shadowing of the quasi-trajectories of the separatrix map.

${\bf 1.}$ In Section 2 we introduce the separatrix map: the construction and formulas (an explicit part plus small error terms). For any point $z$ in a neighbourhood of $\Gamma^+_{\eps} \cup \Gamma^-_{\eps}$ this map is an integer degree of the time-one map of the perturbed system, the degree depends on $z$. The idea is to skip a long noninteresting part of the dynamics, corresponding to a passage through a small neighbourhood of $N_{\eps}$.

We study the separatrix map in the variables $\rho, \zeta, \tau, t$, where $\tau \in [-1,1]$, $t \in \mN$. Up to some small correction terms $\eps \rho=y$ while $\zeta-x$ is a function of $y,u$ and $v$. Dynamical meaning of $\tau$ and $t$ is more complicated. Roughly speaking, $(t+\tau)(z)$ is the time the perturbed flow needs to put the point $z$ to a certain cross section after one homoclinic excursion (i.e. after one passage through a small neighbourhood of $N_{\eps}$).

Important ingredients of the separatrix map (see Theorem \ref{theo:sm_main}) are the Poincar\'e--Melnikov potentials $\hat\Theta^{\pm}$ (the signs $+$ and $-$ distinguish the ``upper'' and ``lower'' homoclinic loops) and the function $\bf H$, obtained as smoothed averaging of $H_1$ at the hyperbolic manifold $N$.

${\bf 2.}$ In Section 3 we define a special class of trajectories (finite pieces)
\begin{equation}
\label{adm_traj0}
(\rho_j,\zeta_j,\tau_j,t_j), \qquad 0 \le j \le m
\end{equation}
of the separatrix map, the class of admissible trajectories. The idea is to keep the dynamics in a $c\eps$-neighbourhood of $\Gamma^+_{\eps} \cup \Gamma^-_{\eps}$, but to avoid coming too close to $\Gamma^{\pm}_{\eps}$ (otherwise the error terms in the separatrix map may start to dominate over the main explicit terms).
For any $j$ this requirement is reduced to the inequality
\begin{equation}
\label{integer_t}
-A'(\eps \rho_j) \log \eps \le t_j \le -A''(\eps \rho_j)\log \eps
\end{equation}
for certain positive functions $A', A''$
and to the assumption that the point $(\zeta_j,\tau_j)$ lies in a neighbourhood of the set of $\JJ_0(\eps\rho_j)$ defined in terms of the Poincar\'e-Melnikov potential $\hat \Theta = \hat \Theta^+$.

By definition any admissible trajectory (\ref{adm_traj0}) has a code, a quasi-trajectory
\begin{equation}
\label{quasi-traj0}
(\overline\rho_j,\overline\zeta_j,\overline\tau_j,\overline t_j), \quad 0 \le j \le m, \quad (\overline \zeta_j,\overline \tau_j) \in \JJ_0(\eps\rho_j),
\end{equation}
close to (\ref{adm_traj0}) in some metric. We present rules for a proper extending of the code i.e., for adding of a new point
$(\overline\rho_{m+1},\overline\zeta_{m+1},\overline\tau_{m+1},\overline t_{m+1})$. According to Lemma \ref{lem:r_attach} the extended code corresponds to some admissible trajectory. In this way, we obtain a piece of trajectory with $m+2$ points.

${\bf 3.}$
Thus, it is sufficient to construct a quasi-trajectory the $\overline \rho$-component of which moves in a small neighbourhood of the curve $\chi$ and to estimate the average velocity of this motion.

The variables $\overline \rho$ and $\overline \zeta$ on a quasi-trajectory satisfy the equations

\begin{eqnarray}
\label{rho01}
\!\!\!\!\!\!\!\!
\overline \rho_{m+1}-\overline\rho_m \!\! &=& \!\! - \hat \Theta_{\zeta}\big(\eps \overline \rho_m,\overline \zeta_m,\overline \tau_m\big)+(\overline\tau_{m+1}-\tau_m-t_{m+1})\bH_{\zeta}\big(\eps\overline\rho_m,\zeta_m\big)+r_m,\\
\label{zeta(rho)01}
\!\!\!\!\!\!\!\!
\overline\zeta_{m+1}-\overline\zeta_m\!\! &=& \!\!\hat\nu(\varepsilon \overline \rho_m)\,t_{m+1}+s_m.
\end{eqnarray}
(compare with (\ref{rho1})--(\ref{zeta(rho)1})). Here
$r_m$, $s_m$ are small error terms and $t_{m+1}$ is an integer number satisfying (\ref{integer_t}).

%$(ii)$ $s_{\beta,a}(\zeta)$ is a piecewise smooth function such that the norm of
%\begin{equation}
%\label{integral0}
%\int_{\mT^n} \hat\Theta_{\zeta}\big(a,\zeta,s_{\beta,a}(\zeta)\big)\, d\zeta-\lambda_{\Theta}\beta
%\end{equation}
%is small for any unit vector $\beta \in \mR^n$ and some positive $\lambda_{\Theta}$ which depends only on $\Theta$.

The main term which pushes $\overline \rho$ in a given direction is $-\hat \Theta_{\zeta}\big(\eps \overline \rho_m,\overline \zeta_m,\overline \tau_m\big)$. The second term in the right-hand side of (\ref{rho01}) vanishes outside small neighbourhoods of essential resonances, i.e. resonances determined by resonance vectors $|h|=O\big(|\log \eps |^{1/({\bf r}-n-1)}\big)$. When $\rho_m$ passes through an $\eps^{1/4}$-neighbourhood of such a resonance, this term is noticeable, but due to the fact that $\zeta$-average of $\bf{H}_{\zeta}$ vanishes, contribution of this term is small on long parts of a trajectory.

The possibility to keep $-\hat \Theta_{\zeta}\big(\eps \overline \rho_m,\overline \zeta_m,\overline \tau_m\big)$ directed along the curve $\chi$ having not very small norm, depends on the possibility to choose the point $(\overline\zeta_m,\overline \tau_m)$ of the code in some special parts of the set $\JJ_0(\eps\overline\rho_m)$.
According to the rules for the extension of the code $(\overline\zeta_m,\overline \tau_m) \in \JJ_0(\eps \overline \rho_m)$ and equations (3.14)--(3.15) should hold.

In other words, with high precision the point
$$
(\overline \xi_{m+1}, \overline \tau_{m+1}) = (\overline \zeta_{m+1} - \hat\nu(\eps \overline \rho_{m+1})\overline \tau_{m+1}, \overline \tau_{m+1}) \in \mT^{n+1}
$$
lies on an interval of a winding through the point $(\overline \xi_m, \overline \tau_m) \in \mT^{n+1}$, where $\overline t_m$ lies on an interval of length $\hat k=|\log \eps|/10$ (see Remark \ref{remark_3_1}) and the frequency vector of the winding is
$\nu(\eps \rho_{m+1})=(\hat\nu(\eps \rho_{m+1}),-1) \in \mR^{n+1}$.

The point $(\overline \xi_{m+1}, \overline \tau_{m+1})$ should lie in the set $J(\eps\overline\rho_m) \subset \mT^{n+1}$, the image of $\JJ(\eps\overline\rho_m)$ under the map $(\zeta,\tau) \mapsto (\xi,\tau)=(\zeta-\hat\nu(\eps \overline \rho_m)\tau,\tau)$.
Properties of the function $\hat\Theta_{\zeta}$ on the set $\JJ(\eps\overline\rho_m)$ are discussed in Section \ref{proofLem}. To have values of $\hat \Theta_{\zeta}$ which push $\overline \rho_m$ in a given direction, the interval of winding should intersect a proper part of $J(\eps \overline \rho_m)$. This can be achieved relatively easily if the point $\eps \overline \rho_{m+1}$ is far from low order resonances. This program was accomplished in \cite{DT}.

${\bf 4.}$
The problem of the passage through a neighbourhood of a strong (low order) resonance was solved in \cite{DavTre16} under additional assumption that $H_1$ is a trigonometric polynomial in $x$. In this paper we consider the general case.
Near a strong resonance
\begin{equation*}
\langle  \nu, k\rangle = 0, \quad \nu = (\hat \nu,-1), \quad k \in \mZ^{n+1}, \quad |k| \le C_{\diamond}
\end{equation*}
the quantities $\langle (\overline \xi_{m+1},\overline \tau_{m+1})-(\overline \xi_m,\overline \tau_m),k\rangle$ are small and we can not hope that the point $(\overline \xi_{m+1},\overline \tau_{m+1})$ can always hit the set $J(\eps \overline \rho_m)$ close to a point where the direction of the vector $\Theta_{\zeta}$ is prescribed. Hence our strategy is as follows. We consider evolution of the variable $\eps \rho$ on a large piece of a trajectory when the slowly varying angle $\langle (\xi,\tau), k\rangle$ mod $1$ makes a full rotation. On such a piece the total contribution of the terms $-\hat\Theta_{\zeta}(\eps \overline \rho_m,\overline \zeta_m,\overline \tau_m)$ in (\ref{rho01}) can be made directed along $\chi$ and providing the required velocity of the diffusion.

We have several obstacles in implementing this program. The first one is
the existence of so-called essential resonances\footnote{Note that if $H_1$ is a trigonometric polynomial in $x$ (\cite{DavTre16}), we can assume that the neighbourhood of the strong resonance we pass does not contain other resonances.}. The corresponding integer resonant vectors $h \in \mZ^{n+1}$ are not very long:
\begin{equation*}
|h|<\Big(\frac{|\log\eps|}{cd\underline\lambda}
    \Big)^{1/({\mathbf r}-n-1)}
\end{equation*}
(see \ref{h_c}), here $c$, $d$, $\underline\lambda$ are some positive constants, that do not depend on $\eps$.
The term $(\overline\tau_{m+1}-\tau_m-t_{m+1})\bH_{\zeta}\big(\eps\overline\rho_m,\zeta_m\big)$ in (\ref{rho01}) may be big as $\eps\overline \rho_m$ is close to an essential resonance and prevents to push the quasi-trajectory in a chosen direction.
First, deforming slightly the curve $\gamma$, we replace it by a polygon line (Lemma \ref{lem:deform}) such that each interval of this line is either sufficiently far from all essential resonances or crosses exactly one $\eps^{1/(7n)}$-neighbourhood of an essential resonance. In this way we avoid small neighbourhoods of essential multiple resonances.

Passage through an $\eps^{1/4}$-neighbourhood of the strong resonance requires an additional argument.
We approximate the discrete system by a second order ODE system (\ref{alpha.ph.2.1}).
It appears to be similar to a pendulum with an additional constant rotational force. We study the dynamics of this system and prove that its solutions approximate well the corresponding solutions of the discrete system. Then we prove that the solutions of the continuous system intersect the strong resonance in Section \ref{Etap2}. Therefore, there exists a quasi-trajectory that intersects the strong resonance.

Below all constants we use except $\overline K = |\log \eps|/10$ are independent of $\eps$. Hence dependence of our estimates on $\eps$ is always explicit.

\section{The separatrix map}
Invariant $(n+1)$-dimensional tori (\ref{N(y)}) are called partially hyperbolic \cite{Gra,Zeh,BT,Tre_book}. The asymptotic manifolds $\Gamma^\pm(y^0)$,
\begin{eqnarray*}
& \Gamma^+(y^0), \Gamma^-(y^0) \subset
  \{(y,x,v,u,t) : y=y^0, H_0(y^0,v,u) = H_0(y^0,0,0) \}, & \\
&   \Gamma^\pm(y^0)
  =  \{y^0\} \times \mT^n \times \gamma^\pm \times \mT.
\end{eqnarray*}
consist of unperturbed solutions which tend to $N(y^0)$ as $t\to\pm\infty$.

Consider the dynamics of the perturbed system in a small neighborhood of the set
$$
    \Gamma
  = \cup_{y\in\overline\DD}
          \big( \Gamma^+(y) \cup \Gamma^-(y) \big).
$$
Let $T_{\eps}$ be the time-one shift.
If $\eps=0$, the shift $T_0$ is an integrable symplectic map, for which  $L(y) = \mbox{pr}(N(y))$ are $n$-dimensional hyperbolic tori and $\Sigma^\pm(y) = \mbox{pr}\,(\Gamma^\pm(y))$ are asymptotic manifolds. Here pr is the projection $(y,x,v,u,t) \mapsto (y,x,v,u)$.

Now we define the separatrix map $\SM_\eps$ corresponding
to $T_\eps$ in the neigbourhood of $\Sigma$,
$$
    \Sigma
  = \cup_{y\in\overline\DD}
      \big( \Sigma^+(y) \cup \Sigma^-(y) \big).
$$

Let $U$ be a small neighborhood of the unperturbed normally hyperbolic manifold $L = \cup_{y\in\overline\DD} L(y)$ and $\bU$ a neighborhood of $\Sigma$. If $\bU$ is sufficiently small, $\bU\setminus U$ breaks into two connected components ${\bU}^+$ and ${\bU}^-$ such that $\Sigma^\pm\subset {\bU}^\pm\cup U$.

Consider a point $z\in{\bU}^+ \cup {\bU}^-$. Let $m_1 = m_1(z)$ be the minimal natural number such that
$T_\eps^{m_1}(z) \not\in {\bU}^+ \cup {\bU}^-$
and let $m_2 = m_2(z)$ be the minimal natural number such that
$m_2 > m_1$ and $T_\eps^{m_2}(z) \in {\bU}^+ \cup {\bU}^-$.
So, for $m=m_1$ the trajectory $T_\eps^m (z)$ leaves the domain
${\bU}^+ \cup {\bU}^-$. For $m=m_2$ the trajectory returns to
${\bU}^+ \cup {\bU}^-$. Denote by $\bG$ the set of points $z$ such that $m_2 < \infty$
and $T_\eps^{m_1}(z),\ldots,T_\eps^{m_2 - 1}(z)\in U$.
Putting
$$
  {\bU}_\eps = ({\bU}^+ \cup {\bU}^-) \cap \bG,
$$
we obtain the map
$$
{\bU}_\eps \ni z\mapsto  \SM_\eps(z) = T_\eps^{m_2(z)}(z) \in {\bU}^+ \cup {\bU}^-.
$$

\subsection{Explicit formulas for the separatrix map}
Consider the Fourier expansion of the function $H_1(y,x,0,0,t)$:
$$
  H_1(y,x,0,0,t) = \sum_{k\in\mZ^n,k_0\in\mZ}
             H_1^{k,k_0}(y) e^{2\pi i(\langle k,x\rangle + k_0 t)}.
$$
Let $\phi : \mR\to [0,1]$ be a $C^\infty$-smooth function such that
$\phi(r) = 0$ for any $|r| \ge 1$, and
$\phi(r) = 1$ for any $|r| \le 1/2$.

We define the following smoothed averages of $H_1$.
\begin{eqnarray}
\label{eq:overlinebH}
\hskip-.5cm
  \overline{\bH}(y,x,t)
  \! &=& \!
       \sum_{(k,k_0)\in \mZ^{n+1}}
             \phi \Big( \frac{ \langle k,\nu(y)\rangle}
                             {\eps^{1/4}}
                  \Big)
             H_1^{k,\,k_0}(y)\, e^{2\pi i(\langle  k,x\rangle  + k_0 t)},\\
\label{eq:bH=}
\hskip-.5cm
  {\bH}(y,x)
  \! &=& \!
            \overline{\bH}(y,x,0).
\end{eqnarray}

\begin{dfn}
For any $f\in C^j(\overline {\cal D}\times \mT^n)$,
$g\in C^0(\overline {\cal D}\times \mT^n)$, and $\delta = \delta(\eps) > 0$ we say that
$f = O^{(\delta)}(g)$ if for any $l',l''\in \{0,1,\ldots\}$,
$l' + l'' := l \le j$
$$
        \Big| \frac{\partial^{l'+l''} f}
                   {\partial y_1^{l'_1}\ldots\partial y_n^{l'_n}
                    \partial x_1^{l''_1}\ldots\partial x_m^{l''_m}}
        \Big|
  < C_l   \delta^{-l'} \, |g|, \qquad
  y\in\overline{\DD},\quad  x\in\mT^n
$$
with $C_l$ independent of $\delta$,
where $l' = l'_1 + \ldots + l'_n$,  $l'' = l''_1 + \ldots + l''_m$.
Here we assume that $f$ can take values in $\mR^s$, where $s$ is
an arbitrary natural number.
\end{dfn}

According to (\ref{eq:bH=})
$\overline{\bH},\bH \in C^{\bf r}(\overline {\cal D}\times \mT^n)$
and $\bH = O^{(\eps^{1/4})}(1)$. More precisely,
\begin{eqnarray}
\label{|H|<}
& |\bH_x| < C_H, \quad
  |\bH_y| < (C_0\eps^{-1/4} + 1)\, C_H, & \\
\nonumber
&   C_H
  = 2\pi\sum_{(k,k_0)\in\mZ^{n+1}}
              (1 + |k|)\, \|H_1^{k,k_0}\|_{C^1}, \quad
    C_0
  = \|H_{0yy}\|_{C^0}. &
\end{eqnarray}

\begin{prop}
\label{stm:bH-quasi}
For any $y\in\overline\DD$, $x\in\mT^n$ and $l\in\mZ$
\begin{equation}
\label{H_periodic}
    \bH(y,x + \hat \nu(y) l)
  = \bH(y,x) + O^{(\eps^{1/4})}(\eps^{1/4} l).
\end{equation}
\end{prop}

In \cite{Tre_sm_gl} the coordinates
$(\rho,\zeta,r,\tau,\sigma)$ on $\bU^+\cup\bU^-$ are constructed such that
\smallskip

1) $dy\wedge dx + dv\wedge du
    = \eps (d\rho\wedge d\zeta + dr\wedge d\tau)$,

2) for some function $f = f(y,u,v,\eps) = O^{(\eps^{1/4})}(1)$ such that $f(y,0,0,0) = 0$,
 \begin{eqnarray*}
   && \eps\rho = y + O^{(\eps^{1/4})}(\eps^{3/4},H_0 - E),\quad
         \zeta = x + f(y,u,v,\eps), \\
   && \eps r   = H_0 - E + O^{(\eps^{1/4})}(\eps^{3/4},H_0 - E),
 \end{eqnarray*}
where $H_0 = H_0(y,u,v)$ and $E = E(y)$,

3) the variable $\tau\in [-1,1]$ is an analogue of time $t$,

4) $\sigma \in \{-1,1\}$ stands for a domain ($1$ corresponds to $\bU^+$ and $-1$ to $\bU^-$),
in which the phase point lies.
\smallskip

We also need another discrete variable
$$
  \thet\in\{-1,1\}, \qquad
  \thet_{m+1} = \mbox{sign} (r_{m+1} - \bH(\eps\rho_{m+1},\zeta_m)).
$$
\begin{theo}
\label{theo:sm_main}
{\rm (\cite{PifTre},\cite{Tre_sm_gl},\cite{Tre_sm_sym})}.
Suppose that Hypotheses ${\bf H_0 1}$--${\bf H_0 2}$ hold,
$\overline K = \frac1{10}|\log\eps|$, $K_0 > 0$
is a (large) constant independent of $\eps$ and satisfying the inequality $\overline K + K_0 < \frac 19 |\log\eps|$.

Then there exist $C^{\bf r}$-smooth functions
$$
  \lambda,\kappa^\pm : \overline {\cal D} \to \mR,\quad
          \Theta^\pm : \overline {\cal D}\times\mT^{n+1}\to\mR
$$
and coordinates $(\rho,\zeta,\tau,t,\sigma,\thet)$ such that
for any trajectory $(\rho_m,\zeta_m,\tau_m,t_m,\sigma_m,\thet_m)$
of the separatrix map, where
\begin{equation}
\label{eq:assume}
      - \overline K - K_0
  \le - \lambda t_{m+1} - \log\eps
  \le - K_0,
\end{equation}
the following equations hold:
\begin{equation}
\label{eq:Routh}
  \begin{array}{l}
\displaystyle{
\rho_{m+1}  = \rho_m - \widehat\Theta_\zeta^{\sigma_m}
                          (\eps\rho_{m+1},\zeta_m,\tau_m)
               + (\tau_{m+1} - \tau_m - t_{m+1})
                        \bH_\zeta(\eps\rho_{m+1},\zeta_m)
               + \widehat{\bO}_2, } \\[1mm]
\displaystyle{
    \zeta_{m+1} = \zeta_m + \hat\nu(\eps \rho_m) t_{m+1}
                    - (\tau_{m+1} - \tau_m - t_{m+1})
                        \bH_\rho(\eps\rho_{m+1},\zeta_m)
                    + \widehat{\bO}_1, } \\[1mm]
\displaystyle{
    \widehat\Theta_\tau^{\sigma_m}(\eps\rho_{m+1},\zeta_m,\tau_m)
     = \frac\lambda\eps
 \Big( \frac{\thet_m}{\kappa^{\sigma_{m-1}}}
              e^{\lambda (\tau_m - \tau_{m-1} - t_m)}
   - \frac{\thet_{m+1}}{\kappa^{\sigma_m}}
              e^{\lambda (\tau_{m+1} - \tau_m - t_{m+1})}
 \Big)
                       + \widehat{\bO}_2, } \\[1mm]
\displaystyle{
    \sigma_{m+1} = \sigma_m \thet_{m+1} . }
  \end{array}
\end{equation}
Here $\lambda,\nu,\kappa^\sigma$ are functions of $\eps\rho_{m+1}$;

\begin{equation}
\widehat\Theta(\eps\rho,\zeta,\tau)
 = \Theta(\eps\rho,\zeta - \hat \nu(\eps\rho)\tau,\tau);
\end{equation}

\begin{equation}
\label{O1O2}
  \widehat{\bO}_1 = O^{(\eps^{-3/4})}(\eps^{7/8}\log^2\eps),\quad
  \widehat{\bO}_2 = O^{(\eps^{-3/4})}(\eps^{1/8}\log^2\eps),\quad
  \mbox{with respect to} \quad \rho_+.
\end{equation}
\end{theo}
The functions $\lambda>0$ and $\kappa^\pm > 0$ are determined by the unperturbed system.
The curves $\widehat\gamma^\pm$ are parameterized by time up to the shift: $t\mapsto t + t_0(y)$. Then
\begin{eqnarray*}
  \Gamma^\sigma(y,\xi,\tau)
               \!\! &=& \!\! (y, \xi + \nu(y)\tau, \gamma^\sigma(y,\tau)).
\end{eqnarray*}
We put
$$
  H^\sigma_*(y,\xi,\tau,t)
    \eqdef
  H_1(\Gamma^\sigma(y,\xi,t),t - \tau)
                      - H_1(y,\xi + \hat\nu(y) t,0,0,t - \tau).
$$
Then $H_*^\sigma (y,\xi,\tau,t)$ exponentially tends to zero
as $t\to\pm\infty$. Define
\begin{equation}
\label{hattheta}
      \Theta^\sigma (y,\xi,\tau)
  \eqdef
     - \int_{-\infty}^{+\infty} H^\sigma_* (y,\xi,\tau,t) \, dt, \quad
     \widehat\Theta^\sigma (y,\zeta,\tau)
  \eqdef     \Theta^\sigma (y,\zeta - \hat \nu(y)\tau,\tau).
\end{equation}
The functions $\Theta^\pm \in C^{\bf r}(\DD\times \mT^n \times D \times \mT)$ are the Poincar\'e--Melnikov integrals.
The genericity of $H_1$ is equivalent to the genericity of $\Theta^\pm$ in the following sense: a typical set of
Poincar\'e--Melnikov integrals corresponds to some typical set of perturbations. Thus it is sufficient to prove the
existence of the diffusion for a typical set of functions $\Theta^\pm$.

\section{Construction of a trajectory}
\label{sec:symbolic}

In this section we briefly recall results of \cite{DT} (see also \cite{Tre_sm_sym}) about construction of trajectories of the separatrix map. We define a quasi-trajectory (code) and present the rules which let us to extend the code. Lemma \ref{lem:r_attach} states that any extended code generates an extended trajectory.

For simplicity we put $\sigma\equiv 1$ and hence $\thet \equiv 1$.
Together with variables $\rho, \zeta, \tau$
we use variables $\eta,\xi,\tau$ defined by the map
$$
  \pi:\mR^n\times\mT^n\times\mR \to \mR\times\mT^{n+1}, \quad
  \pi(\rho,\zeta,\tau) = (\eta,\xi,\tau)
                       = (\eps\rho,\zeta - \hat\nu(\eps\rho)\tau,\tau),
$$
$$
      \widehat\Theta^\sigma (y,\zeta,\tau)
  =     \Theta^\sigma (y,\zeta - \hat \nu(y)\tau,\tau).
$$
The Poincar\'e--Melnikov integral
$\widehat\Theta = \widehat\Theta^+$ in the new variables turns to $\Theta = \Theta^+$ which is periodic in $\xi$
and $\tau$.
We put
\begin{equation}
\label{nu,partial}
   \partial
 = - \langle\hat\nu,\partial/\partial\xi\rangle + \partial/\partial\tau.
\end{equation}

\begin{dfn}
Consider the sets
\begin{eqnarray*}
     J_0
 &=& \big\{ (\eta,\xi,\tau)\in\overline\DD\times\mT^{n+1} :
            \, \partial\Theta(\eta,\xi,\tau) = 0,
            \, \partial^2\Theta(\eta,\xi,\tau) \ne 0
     \big\}, \\
     \JJ_0
 &=& \pi^{-1} (J_0)
     \cap \big\{ (\rho,\zeta,\tau) : -1 < \tau < 1 \big\}
 \subset \frac1\eps\overline\DD\times\mT^n \times (-1,1).
\end{eqnarray*}
\end{dfn}
The set $\JJ_0$ can be also defined as:
\begin{equation}
\label{eq:calN0}
   \JJ_0
 = \Big\{ (\rho,\zeta,\tau)
            \in \frac1\eps\overline\DD\times\mT^n \times (-1,1) :
                  \widehat\Theta_\tau(\eps\rho,\zeta,\tau) = 0,
              \widehat\Theta_{\tau\tau}(\eps\rho,\zeta,\tau) \ne 0
   \Big\}.
\end{equation}
Consider the equation
\begin{equation}
\label{eq:Theta=z}
  \widehat\Theta_\tau(\eps\rho,\zeta,\tau) = z, \qquad
  z\in\mR.
\end{equation}
It can be solved with respect to $\tau$ for small $|z|$ near any point $(\rho_0,\zeta_0,\tau_0)\in \JJ_0$. The solution is a smooth function $\Psi^{\rho_0,\zeta_0,\tau_0}(\eps\rho,\zeta,z)$ with values in $(-2,2)$.

\begin{dfn}
\begin{eqnarray}
\label{JJc'c''}
    \JJ_{c',c''}
 \!\! &=& \!\!
     \big\{
         (\rho_0,\zeta_0,\tau_0)\in \JJ_0
       : \Psi = \Psi^{\rho_0,\zeta_0,\tau_0}(\eps\rho,\zeta,z)\quad
         \mbox{is smooth for}\\
 && \;
        \quad
        \eps^{3/4} |\rho - \rho_0| < c',\,
        |\zeta - \zeta_0| < c',\,
        |z| < c',  \mbox{ where }\\
 && \;
        |\Psi| < 2 ,\;
  |\Psi_\rho| < \eps^{3/4}/c'' ,\;
 |\Psi_\zeta| < 1/c'' ,\;
     |\Psi_z| < 1/c''
     \big\} , \\
      J_{c',c''}
 \!\! &=& \!\!
      \big\{ (\eta,\xi,\tau)\in J_0
       : \pi^{-1}(\eta,\xi,\tau) \cap \JJ_{c',c''} \ne \emptyset
      \big\} .
\end{eqnarray}
\end{dfn}
It is obvious, that
$$\displaystyle{
  \bigcup_{c'>0,c''>0} \JJ_{c',c''} = \JJ_0, \quad
  \bigcup_{c'>0,c''>0} J_{c',c''} = J_0 .}
$$
Let $c\,',c\,''$ be a fixed sufficiently small constants.
Denote
$$
  \OO = (\Omega_0,\ldots,\Omega_m),\quad
  \Omega_j = (\rho_j,\zeta_j,\tau_j,t_j),\qquad
  0 \le j\le m.
$$
\begin{dfn}
Let $C$ be a big constant and $b$ satisfies
$$
  0 < b < \min\{1/3,c'/2\}.
$$
We define
$$
  b_\rho = \frac{b^5}{60 C^3} e^{K_0},\quad
  b_\tau = \frac{b^4}{3 C^2} e^{K_0},\quad
  b_\zeta = \frac{b^5}{48 C^3} e^{K_0}.
$$
Here we assume that $K_0 = K_0(C,b)$ is chosen so that $b_\rho$, $b_\tau$, and $b_\zeta$ are large.
\end{dfn}

\begin{dfn}
\label{dist}
For two points $\Omega', \Omega''$ define
$$
  \mbox{d}(\Omega',\Omega'')
  = \left\{ \begin{array}{l}
             \; +\infty, \quad \mbox{if}\quad  t'\ne t'', \\
             \max \big\{ b_\rho  |\rho' - \rho''| ,\,
                     \eps^{-3/4} b_\zeta |\zeta' - \zeta''| ,\,
                     b_\tau  |\tau' - \tau''|    \big\} \quad \mbox{otherwise}.
            \end{array}
    \right.
$$
\end{dfn}
Here $| \cdot|$ is a standard metric on the corresponding space.

\begin{dfn}
\label{dfn:one}
We say that
$\overline\OO = (\overline\Omega_0,\ldots,\overline\Omega_m),\;$
$m\ge 0$ is a quasi-trajectory if
\begin{eqnarray}
\label{qt1}
 && (\overline\rho_j,\overline\zeta_j,\overline\tau_j)
           \in \JJ_{c',c''},
    \qquad\qquad\qquad\qquad    0 \le j\le m, \\
\label{qt2}
 && \overline t_j\in\mN , \; |\overline\tau_j| < 1,
   \qquad\qquad\qquad\qquad\quad   0 \le j\le m, \\
\label{qt3}
 &&\!\!\!\!\!\!
    K_0  \le \lambda(\eps\rho_j)\overline t_j + \log\eps
         \le K_0 + \overline K,
   \qquad\quad  0 < j \le m .
\end{eqnarray}
About $K_0$, $\overline K$ see Theorem \ref{theo:sm_main}.
\end{dfn}

\begin{dfn}
\label{dfn:two}
We call a trajectory $\OO$ admissible if a quasi-trajectory (a code)
$\overline\OO$ exists such that
\begin{eqnarray}
\label{admiss:dist}
&&\!\!\!\!
     \mbox{d}(\Omega_j,\overline\Omega_j)
  < b\, (2-b^{1+m-j}), \quad
    0 \le j\le m ,  \\
\label{admiss:incl}
&&    (\rho_0,\zeta_0,\tau_0)
  \in \JJ_{c',c''}, \quad
      (\rho_m,\zeta_m,\tau_m)
  \in \JJ_{c',c''}.
\end{eqnarray}
\end{dfn}

Since $2b < c'$, it follows that inclusions (\ref{admiss:incl}) imply the equations

$$
    \tau_0
  = \Psi^{\overline\rho_0,\overline\zeta_0,\overline\tau_0}
                  (\eps\rho_0,\zeta_0,0), \quad
    \tau_m
  = \Psi^{\overline\rho_m,\overline\zeta_m,\overline\tau_m}
                  (\eps\rho_m,\zeta_m,0).
$$

\begin{dfn}
\label{dfn:three}
Let $\OO$ be an admissible trajectory with code $\overline\OO$.
We say that a quasi-trajectory $(\Omega,\Omega_+)$
is compatible with $\OO$ if
\begin{eqnarray}
\label{eq:fromtheright-}
&&\!\!\!\!\!\!\!\!\!\!
    \Omega = \overline\Omega_m,\\
\label{eq:fromtheright}
&&\!\!\!\!\!\!\!\!\!\!\!\!
  \big| \rho_+
       - \rho_m + \widehat\Theta_\zeta
                               (\eps\rho_+,\zeta_m,\tau_m)
                - (\tau_+ - \tau_m - t_+) \bH_\zeta(\eps\rho_+,\zeta_m)
  \big| < \frac{b^2}{2 b_\rho} ,\\
\label{eq:fromtheright+}
&&\!\!\!\!\!\!\!\!\!\!
  \big| \zeta_+
      - \zeta_m - \hat\nu(\eps\rho_+) t_+
                + (\tau_+ - \tau_m - t_+) \bH_\rho(\eps\rho_+,\zeta_m)
  \big| < \frac{b^2\eps^{3/4}}{2 b_\zeta} .
\end{eqnarray}
\end{dfn}

\begin{rem}
\label{remark_3_1}
Equation (\ref{eq:fromtheright+}) means that the points
$$
(\xi_m,\tau_m)=(\zeta_m-\hat\nu(\eps\rho_m)\tau_m,\tau_m) \in \mT^{n+1} \quad \text{and} \quad
(\xi_+,\tau_+)=(\zeta_+-\hat \nu(\eps \rho_+)\tau_+,\tau_+) \in \mT^{n+1}
$$
satisfy the equation
$$
(\xi_+,\tau_+)=(\xi_m,\tau_m)+ \hat\nu(\eps \rho_m)(\tau_m-\tau_++t_+)+O(\eps^{3/4}) \mod 1.
$$
\end{rem}

We need the following lemma (see \cite{DT,Tre_sm_sym}) to construct a diffusion trajectory.
\begin{lem}
\label{lem:r_attach}
Let $\OO^0 = (\Omega_0^0,\ldots,\Omega_m^0)$ be an admissible trajectory with code $\overline\OO$, and let the quasi-trajectory
$(\overline\Omega_m,\Omega_+)$ be compatible with $\OO^0$.
Then there exists an admissible trajectory
$\widehat\OO = (\widehat\Omega_0,\ldots,\widehat\Omega_{m+1})$
with code
$$
   (\overline\Omega_0,\ldots,\overline\Omega_{m+1}),\qquad
   \overline\Omega_{m+1} = \Omega_+.
$$
\end{lem}

Thus, to extend an admissible trajectory it is sufficient to find a compatible quasi-trajectory.
\section{Essential resonances}
\label{essentialresonance}
In this section we define essential resonances and present a lemma which shows that the curve $\chi$ can be
slightly deformed so that the new curve $\varkappa$ avoids $\eps^{1/(7n)}$-neighbourhoods of multiple essential
resonances.

\begin{dfn}
\label{dfn:essential}
We call a point $\eta\in\DD$ $d$-essential, if
\begin{equation}
\label{H:d}
  \frac{|\log\eps|}{\lambda(\eta)}
   \max_{\zeta\in\mT^n} |\bH_\zeta(\eta,\zeta)| > d.
\end{equation}
\end{dfn}
Let $\EE_d \subset \DD$ be the set of $d$-essential points. In the next section we specify the value of $d$. We will take $d=\lambda_{\Theta}$, where $\lambda_{\Theta}$ is defined in Lemma \ref{lem:typ1}. Essential points lie near resonances of not very big order.
More precisely, in \cite{DT} the following proposition is proven:
\begin{prop}
\label{prop:notverylarge}
Let $\eta\in\EE_d$. Then there exists a constant $c>0$, which
depends only on $C_H = \|H_1\|_{C^{\bf r}}$, $\br$ and $n$ such
that\footnote{If $\br$ equals $\infty$ or $\omega$ we can take in
(\ref{h_c}) instead of $\br$ any integer number, greater than
$n+3$}
\begin{equation*}
\label{notverylarge}
   |\langle \nu(\eta),h\rangle|\le \eps^{1/4} \quad
   \mbox{for some } h\in\mZ^{n+1} \setminus \{0\}, \quad
   |h| < h_c(d), \quad \mbox{where}
\end{equation*}
\begin{equation}
\label{h_c}
    h_c(d)
  = \Big(\frac{|\log\eps|}{cd\underline\lambda}
    \Big)^{1/({\mathbf r}-n-1)}, \quad \underline\lambda=\min_{\DD} \lambda(\eta).
\end{equation}
\end{prop}
%Here $\underline\lambda=\min_{\DD} \lambda(\eta)$.
We call resonances $S^h_0$ with $|h|<h_c(d)$ essential. In particular, all strong resonances are essential.
We put
\begin{equation}
\label{Ess_and_Nonres}
    Z \eqdef \{h=(\hat h,h_0)\in\mZ^{n+1}\setminus\{0\},|h|<|\log \eps|\}, \quad S_{\NN}=S^k_{\delta}\setminus \bigcup_{h \in Z} S^h_{\eps^{1/(6n)}}.
\end{equation}
The domain $S_{\NN}$ in the space of slow variables can be regarded as "nonresonant".
Since $\br>n+2$, then for sufficiently small $\eps$ all $h$ from Proposition \ref{notverylarge} belong to $Z$.

Let $\pi\subset S^k_{\delta}$ be a curve and let $\nu'(\eta) = \partial \nu/\partial \eta$ be the Jacobi matrix. Consider the neighbourhood:
\begin{eqnarray*}
O(\pi,\epsilon)= \{\eta\in S^k_{\delta}: \dist(\eta,\pi)\le\epsilon\}.
\end{eqnarray*}
Suppose that  $h \in Z$ and $\pi$ crosses an essential resonance $S^h_0$ at the point $\eta_0$.
Consider a unit vector $\beta \in \mR^n$ such that
$$
\langle \hat \nu'(\eta_0)\hat h,\beta\rangle>|\hat\nu'(\eta_0)\hat h|/2.
$$

\begin{lem}
\label{lem:deform}
There exists a polygon line $\varkappa\subset S^k_{\delta}$ with endpoints  $\varkappa_0$, $\varkappa_1$ such that

(1) $|\varkappa| < 2 |\chi|$,

(2) $\varkappa \subset O(\chi,\eps^{1/(8n)})$ and $|\varkappa_i-\chi_i|<\eps^{1/(8n)}$, $i=0,1$,

(3) each segment $\pi$ of the polygon line $\varkappa$ is one of the following types:

\quad (a) $O(\pi,\eps^{1/(7n)}) \subset S_{\NN}$, $|\pi|>\eps^{1/(8n)}$, where $|\pi|$ is a length of the segment $\pi$;

\quad (b) $\pi$ intersects an essential resonance $S^h_0$, $h \in Z$ at $\eta_0$. Moreover, $\pi \parallel \beta$ and
$O(\pi,\eps^{1/(7n)})\cap S^{j}_{\eps^{1/(7n)}} = \emptyset$ for all $j \in Z$ such that $j \nparallel h$.
\end{lem}
This lemma is proven in \cite{DT} (Lemma 7.1).
\section{The vicinity of the strong resonance}
\label{5res}
In this section we discuss the behaviour of a typical function $\Theta$ in the vicinity of the strong resonance $S^k_{\delta}$. Recall that $\delta=O(|\log^{-1} \eps|)$.
There exists a constant $C(k)$, which does not depend on $\eps$, such that
$$
S^k_{\delta} \subset \{\eta \colon \inf_{\eta_{\star} \in S^k_0} |\eta-\eta_{\star}|\leq C(k)\delta\}.
$$
By Lemma \ref{lem:deform} we approximate the curve $\chi$ by a polygonal chain $\varkappa\subset S^k_{\delta}$ that does not lie too close to multiple essential resonances.
Since $\chi$ does not contain multiple strong resonances, its small\footnote{but independent of $\eps$} neighbourhood does not contain any strong resonances except $S^k_0$.
Let $\eta_*=\eps\rho_*$ be the point of intersection of the polygonal chain $\varkappa$ and $S^k_0$, and let $\beta \in \mR^n$ be an arbitrary unit vector.
Consider any point $p$ in the ball:
\begin{equation}
\label{p}
p \in \{\eta\colon |\eta-\eta_*|\le C(k)\delta\}.
\end{equation}
Let $\{x\}$ be the fractional part of $x$.
\begin{lem}
\label{lem:typ1}
For each $\Theta$ from a typical subset of $C^{\bf r}(\DD\times \mT^{n+1})$ and each $\epsilon>0$, for sufficiently small positive $\eps,c\,',c\,''$, there exists a function $s_{\beta,p}=s_{\beta,p}(\zeta) \colon \mT^n \to \big[0,(\overline K/\lambda)^{1/2}\big]$ such that

1) $\big(p,\zeta-\hat\nu(p)s_{\beta,p}(\zeta), \{s_{\beta,p}(\zeta)\}\big) \in J_{2c',2c''}$;

2) the set of break points of $s_{\beta,p}$ lies in a finite union of smooth compact submanifolds of codimension $1$ in the torus $\mT^n$

3) $\Big|\int_{\mT^n} \hat\Theta_{\zeta}\big(p,\zeta,s_{\beta,p}(\zeta)\big)\, d\zeta-\lambda_{\Theta}\beta\Big|<\epsilon$.

Here $\lambda_{\Theta}$ is a positive number which does not depend on $\epsilon,\eps,c\,',c\,''$.
\end{lem}
The proof of this lemma is contained in Section \ref{proofLem}.
We use Lemma \ref{lem:typ1} to construct a finite piece of the trajectory (in fact, a finite sequence of codes), which pushes the $\rho$-component in the direction of $\beta$.

Now we can define the set of essential points $\EE$. We call the point essential if it is $d$-essential for $d= \lambda_{\Theta}$.
Since $-\hat \Theta^{\sigma_m}_{\zeta}$ in (\ref{eq:Routh}) is an important part of $\rho_{m+1}-\rho_m$, the functions $s_{\beta,p}$ will be useful to control the increment $\langle\rho_{m+1}-\rho_m,\beta\rangle$, which will characterise the motion of the trajectory along $\varkappa$.

We put
\begin{equation}
\label{F}
\mathcal F_{\beta,p}(\zeta)=\hat \Theta_{\zeta}\big(p,\zeta,s_{\beta,p}(\zeta)\big), \quad \calA(p)=\{\zeta \in \mT^n \colon \mathcal F_{\beta,p}(\cdot)\notin C^0(\zeta)\}.
\end{equation}
In other words, $\calA(p)$ is a set of break points of $\mathcal F_{\beta,p}(\cdot)$.

Bellow we assume that the admissible trajectory $\OO$ has already entered
$S^k_{\delta}$.
The next lemma shows how we construct a compatible quasi-trajectory. Then by Lemma \ref{lem:r_attach} we obtain an admissible trajectory that corresponds to the longer code.

\begin{lem}
\label{lem-main}
Consider an admissible trajectory $\OO$ with a code $\overline \OO$
and a point $p$ from (\ref{p}). Let $\overline \rho_m$ satisfies
the inequality $|p - \eps \overline \rho_m|< \eps^{1/(10n)}$ and
let $\overline t$ be an integer number such that
\begin{equation}
\label{overlinetm}
\big(K_0-\log\varepsilon\big) /\lambda + \big(\overline K/\lambda\big)^{1/2}\leq \overline t\leq \big(K_0 - \log\varepsilon + \overline K \big)/\lambda, \quad  \lambda = \lambda(p).
\end{equation}
Then for every unit vector $\beta \in \mR^n$ there exists a compatible quasi-trajectory $(\overline \Omega_m, \overline \Omega_{m+1})$ such that
\begin{equation}
\label{eq1}
\hat \Theta_{\zeta}(\varepsilon \overline\rho_{m+1},\overline\zeta_{m+1},\overline\tau_{m+1})=\mathcal F_{\beta,p}\big(\zeta_m+\hat\nu(p)\overline t\big)+\alpha_{\beta},
\end{equation}
\begin{equation*}
\quad |\alpha_{\beta}|=O\big(\eps^{1/(10n)}|\log^{-1} \varepsilon |\big).
\end{equation*}
\end{lem}
\textbf{Proof.}
For brevity we denote $\Omega_+=\overline\Omega_{m+1}$.
By Lemma \ref{lem:typ1} we define
\begin{equation}
\label{overlinet1}
t \in \Big [-\big(\overline K/\lambda\big)^{1/2},0\Big], \quad \tau_m - t = s_{\beta,p}(\zeta_m+\hat \nu(\eps\rho_m)\overline t).
\end{equation}
Then $\big(\varepsilon \rho_m,\zeta_m+\hat\nu(\varepsilon\rho_m)( t+\overline t-\tau_m),\tau_m- t-\overline t\big) \in J_{2c',2c''}$.
We define $t_+,\tau_{\beta},\zeta_{\beta}$ by formulas
$$
t_+-\tau_{\beta}=t+\overline t-\tau_m, \quad t_+ \in \mZ, \quad \tau_{\beta} \in (-1/2,1/2], \quad \zeta_{\beta}=\zeta_m+\hat\nu(\varepsilon \rho_m)t_+.
$$
Then
$t_+ \in \Big [\big(K_0-\log \varepsilon\big)/\lambda,\big(K_0 - \log \varepsilon + \overline K\big)/\lambda \Big ] $
and $\big(\varepsilon \rho_m,\zeta_{\beta}-\hat\nu(\varepsilon \rho_m)\,\tau_{\beta},\tau_{\beta}\big) \in J_{2c',2c''}$.
Define $\rho_+$, $\zeta_+$ by the equations
\begin{align}
\label{newrho}
\rho_+&=\rho_m-\hat \Theta_{\zeta}(\varepsilon \rho_m,\zeta_m,\tau_m)+(\tau_{\beta}-\tau_m-t_+)\bH_{\zeta}\big(\varepsilon \rho_m,\zeta_m\big),\\
\label{newzeta}
\zeta_+&=\zeta_m+\hat\nu(\varepsilon \rho_m)t_+-(\tau_{\beta}-\tau_m-t_+)
\bH_{\rho}\big(\varepsilon \rho_m,\zeta_m\big).
\end{align}
Find $\tau_+$ from the equation $\hat \Theta_{\tau}(\varepsilon\rho_+,\zeta_+,\tau_+)=0$, then $\tau_+=\Psi^{\rho_m,\zeta_{\beta},\tau_{\beta}}(\varepsilon\rho_+,\zeta_+,0)$.
Since $\tau_{\beta}=\Psi^{\rho_m,\zeta_{\beta},\tau_{\beta}}(\varepsilon\rho_m,\zeta_{\beta},0)$ and $|\zeta_+-\zeta_{\beta}|=O(\varepsilon^{3/4}|\log \varepsilon|)$, then
$|\tau_+-\tau_{\beta}|=O(\varepsilon^{3/4}|\log \varepsilon|)$ and $|\tau_+|<1$.
Thus,
\begin{align*}
\left|\hat\Theta_{\zeta}(\varepsilon\rho_+,\zeta_+,\tau_+)-\mathcal F_{\beta,p}\big(\zeta_m+\hat\nu(p)\overline t_{m+1}\big)\right| \le \left|\hat\Theta_{\zeta}(\varepsilon\rho_+,\zeta_+,\tau_+)-\hat \Theta_\zeta(\varepsilon \rho_m,\zeta_{\beta},\tau_{\beta})\right|\\
+\left|\hat \Theta_\zeta(\varepsilon \rho_m,\zeta_{\beta},\tau_{\beta})-\mathcal F_{\beta,p}\big(\zeta_m+\hat\nu(p)\overline t_{m+1}\big)\right|=O\big(\eps^{1/(10n)}|\log^{-1}\varepsilon|\big).
\end{align*}
Therefore, inequalities  (\ref{eq:fromtheright+}) and (\ref{eq:fromtheright}) hold,
$(\rho_+,\zeta_+,\tau_+) \in \JJ_{2c',2c''}$ and $t_+ \in \mZ$, $K_0 \leq \lambda t_++\log \varepsilon \leq K_0+\overline K$.
Thus $(\overline\Omega_m,\overline\Omega_+)$ is a compatible quasi-trajectory.
\qed

\section{Proof of Lemma \ref{lem:typ1}}
\label{proofLem}
Recall, that $k=(\hat k,k_0) \in \mZ^{n+1}\setminus\{0\}$ is an integer vector, which corresponds to the strong resonance, $\eta_*=\eps \rho_*$ is the point of intersection of the polygonal curve $\varkappa$ and the manifold $S^k_0$. For brevity, we use the notation $\hat \nu_* =\hat \nu(\eta_*)$, $ \hat \nu_* = \hat \nu(\eta_*)$. We define the set
\begin{equation}
\JJ_{2c',2c''}(\eta_*)=\{(\zeta,\tau) \colon (\eta_*,\zeta,\tau) \in \JJ_{2c',2c''}\}.
\end{equation}
The torus $\mT^{n+1}$ is foliated by the $n$-tori
\begin{equation}
\label{mT_theta}
   \mT^n_{\thet}=\{ (\zeta,\tau)\in\mT^{n+1} \colon
       \langle \hat k,\zeta\rangle + k_0\tau = \thet = \mbox{const}\}.
\end{equation}
{\bf 1.}  First we show that for sufficiently small positive $c\,',c\,''$ and every point $\zeta \in \mT^n$ there exists a number $\tau \in\Big[0,\big(\overline K/\lambda\big)^{1/2}\Big]$ such that $(\zeta-\hat \nu_*\tau,\tau) \in \JJ_{2c',2c''}(\eta_*)$.

${\bf H_1 3}$. {\it
For every $\thet \in \mT^1$ there exists a point $\big(\zeta_0(\thet),\tau_0(\thet)\big) \in \mT^n_{\thet}$ such that
\begin{equation}
   \partial\Theta(\eta_*,\zeta_0(\thet),\tau_0(\thet))=0, \quad  \partial^2\Theta(\eta_*,\zeta_0(\thet),\tau_0(\thet)) \ne 0.
\end{equation}}
The operator $\partial$ is defined in (\ref{nu,partial}).

\begin{rem}
1) ${\bf H_13}$ holds for a $C^{\bf r}$-open dense set of functions $\Theta$.

2) ${\bf H_13}$ holds for a subset of full measure in any typical in the sense $1)$ family of functions $\Theta$, which depends on one or several parameters.
\end{rem}
We assume that $\Theta$ satisfies this hypothesis.
\begin{lem}
\label{prop:THgeneric2}
For sufficiently small positive $c\,'$ and $c\,''$ the set
$$
\bigcup_{t \in [0,1]} g^t_{ \nu_*} \big( \JJ_{2c',2c''}(\eta_*)\big)  \cap \mT^n_{\thet}
$$
contains an $n$-dimensional ball of radius $R_{\Theta}$. The radius $R_{\Theta}=R_{\Theta}(\thet)$ depends on  $c\,',c\,'',\|H_0\|_{C^1}$, $\|\Theta\|_{C^2}$, $\thet$, and $n$, but independent of $\eps$.
\end{lem}
The proof of this lemma contained in \cite{DT} (Lemma 6.1).

\begin{dfn}
\label{def:cl-va}
Let $d,K$ be positive numbers. We call the point $\eta\in\overline{\DD}$ $(d,K)$-vague, if for every $n$-dimensional ball $B_d\subset\mT^n_{\thet}$ of the radius $d$
\begin{equation}
\label{Ug(B)=T}
    \bigcup_{0\le t\le K} g^t_{ \nu(\eta)}(B_d)
 \neq \mT^n_{\thet}.
\end{equation}
\end{dfn}

The following proposition says that vague points lie near multiple resonances of small orders.
\begin{prop}
\label{k,epsilon}
Let the point $\eta\in S^k_0$ be $(d,K)$-vague and
$\epsilon  < 1 / (2\pi)$, where
\begin{equation*}
    \epsilon
  = \frac{4n\log (2l_{d,n} + 1)}{\pi K}, \quad
    l_{d,n}
  = \frac{3\cdot 2^n}{\pi d_{\star}^n \sin (\pi d_{\star}/2)} + 1, \quad
  d_{\star}=\frac{d}{(n+1)^{n/2}|k|^n_{\infty}}
\end{equation*}
Then there exists $l\in\mZ^{n+1}$, $l \nparallel k$, $| l|_\infty \le 2\sqrt{n+1}\,l_{d,n}|k|_2$ such that $\eta \in S^l_\epsilon$.
\end{prop}
We prove the proposition in Appendix \ref{section:vague}.

Now we apply Proposition \ref{k,epsilon} for $d=R_{\Theta}$, $K=\big(\overline K/\lambda\big)^{1/2}-1$. Since the point $\eta_*$ is not a multiple resonance, then it is not $(d,K)$-vague for small $\eps$. Therefore,

\begin{equation}
    \bigcup_{0\le t\le K} g^t_{\nu_*}(B_d) = \mT^n_{\thet}.
\end{equation}
Hence the part of winding (with frequency vector $ \nu(\eta_*)$) $(\zeta(t),\tau(t))$ of every point $(\zeta,\tau) \in \mT^{n+1}$ intersects the hypersurface $\JJ_{2c',2c''}$.

For all $p \in B(\eta_*,C(k)\delta)$ we have
\begin{equation}
\label{g-flow}
\big|g^t_{ \nu(p)}(\xi,t)-g^t_{\nu_*}(\xi,t)\big|< C |\log \eps|^{-1/2}.
\end{equation}
Hence for sufficiently small $\eps>0$ for all $\zeta \in \mT$ there exists $\tau \in\Big[0,\big(\overline K/\lambda\big)^{1/2}\Big]$ such that $(\zeta-\nu(p)\tau,\tau) \in \JJ_{2c',2c''}(p)$.

\medskip
{\bf 2.} Now we construct the functions $s_{\beta,p}(\cdot)$.
Let be $\WW$ the set of piecewise smooth maps
\begin{equation}
\label{pw_smooth}
 w : \mT^1\to \JJ_0, \quad
  \thet\mapsto (\zeta,\tau) = w(\thet) \in\mT_{\thet}^n \cap \JJ_0.
\end{equation}
The words ``piecewisely smooth'' mean that there exists a finite collection of intervals $I_j\subset\mT^1$ such that $\mT^1=\cup\overline I_j$ and $w|_{I_j}$ are smooth. This function  is, in general, discontinuous at a finite number of points (end-points of the intervals $I_j$).

Given $w_{\beta} \in \WW$ we define
\begin{equation}
\label{G}
G(\thet)=\Theta_{\zeta}(\eta_*,\zeta(\thet),\tau(\thet)), \quad (\zeta(\thet),\tau(\thet))=w_{\beta}(\thet).
\end{equation}

${\bf H_1 4}$. {\it
There exists a number
$\omega_{\Theta} > 0$ such that for all $\beta\in\mR^n$, $|\beta|=1$ and for all $\omega \in [0,\omega_{\Theta}]$ the following equation holds
$$
\int_{\mT^1} G(\thet)\, d\thet = \omega\beta \quad \text{for some} \quad w_{\beta}\in\WW.
$$
}

\begin{lem}
\label{lem:FF}
${\bf H_1 4}$ holds for generic functions $\Theta$.
\end{lem}
We prove Lemma \ref{lem:FF} in Appendix \ref{subsec:stres}.

Bellow we assume that the function $\Theta(\eta_*,\cdot,\cdot) \colon \mT^{n+1} \to \mR$ satisfies ${\bf H_1 4}$.
Then there exists a piecewise smooth map $w_{\beta}\in\WW$ such that for some $\omega_{\Theta}>0$
\begin{equation*}
\int_{\mT^1} G(\thet)\, d\thet = \omega_{\Theta}\beta.
\end{equation*}
Let us denote the set of discontinuity of $G$ by $\Phi$:
\begin{equation}
\label{Phi}
\Phi=\{\thet=\langle \zeta,\hat k\rangle \in \mT^1 \colon G \notin C^0(\thet)\}.
\end{equation}
For sufficiently small $\epsilon_0>0$ the set
\begin{equation}
\bigcup_{t \in [0,1]} g_{ \nu_*}^t \big( \JJ_{2c',2c''}(\eta_*)\big)  \cap \mT^n_{\thet}
\end{equation}
contains a ball $B(w_{\beta}(\thet),\epsilon_0)$ of radius $\epsilon_0$ and center $w_{\beta}(\thet)$ for all points of continuity of the function $w_{\beta}(\thet)$. Then there exists $\eps_0>0$ such that for all $\eps \in (0,\eps_0)$
\begin{equation}
    \bigcup_{0\le t\le \big(\overline K/\lambda\big)^{1/2}} g^t_{ \nu_*}B\big(w_{\beta}(\thet),\epsilon_0\big) = \mT^n_{\thet}.
\end{equation}
Thus, the winding $(\zeta(t),t)$, $0 \le t \le \big(\overline K/\lambda\big)^{1/2}$ with frequency vector $ \nu(\eta_*)$ with any initial point  $(\zeta,0) \in \mT^n_{\thet}$ intersects the ball $B\big(w_{\beta}(\ph),\epsilon_0\big)$ for all $\thet \in I_j$.

Consider an arbitrary point $p\in B(\eta_*,C(k)\delta)$. By (\ref{g-flow}) for sufficiently small $\eps$ the winding with frequency vector $ \nu(p)$ with any initial  $(\zeta,0) \in \mT^n_{\thet}$ also intersects the ball $B\big(w_{\beta}(\thet),\epsilon_0\big)$ at the time moment $\tau_p(\zeta) \le \big(\overline K/\lambda\big)^{1/2}$ for all $\thet \in I_j$. Consider a function $s_{\beta,p}(\zeta)\eqdef \tau_p(\zeta)$. It is defined on the set
$$
\mathcal C=\big\{\zeta \colon (\zeta,0) \in \mT^n_{\thet}, \quad \thet \notin \cup_j I_j\big\}.
$$
The set $\mT^n \setminus \mathcal C$ is a finite union of compact submanifolds of codimension $1$.
For $\zeta \in \mT^n\setminus \mathcal C$ we put $s_{\beta,p}(\zeta)=0$.

The function $s_{\beta,p}$ may be discontinuous in the following cases:
either if $\thet \notin \cup_jI_j$, or if $\tau_p(\zeta) = 0$ or if $\tau_p(\zeta)=\big(\overline K/\lambda\big)^{1/2}$.
Since
\begin{equation}
\label{g(ph)-h(ph)}
\Big|G(\langle\zeta,\hat k\rangle)-\hat\Theta_{\zeta}\big(p,\zeta,s_{\beta,p}(\zeta)\big)\Big|<c_{\star}\epsilon_0
\end{equation}
for $\zeta \in \mathcal C$, then
\begin{equation}
\Bigg|\int_{\mT^n} \hat\Theta_{\zeta}\big(p,\zeta,s_{\beta,p}(\zeta)\big)\, d\zeta-\omega_{\Theta}\vol(\mT^n_{\ph,0})\beta\Bigg|<c_{\star}\epsilon_0.
\end{equation}
Here $c_{\star}=2\big\|\Theta(\eta_*,\cdot,\cdot)\big\|_{C^2(\mT^{n+1})}$,
and $\vol(\mT^n_{\ph,0})\ge 1$ is the $(n-1)$-volume of the torus
$$
\mT^n_{\thet,0}=\{ \zeta\in\mT^n \colon
       \langle \hat k,\zeta\rangle = \thet = \mbox{const}\}.
$$
Then for $\epsilon_0 <\epsilon/c_{\star}$ we obtain Lemma \ref{lem:typ1}.

\section{Passage through $S_{\delta}^k$}
In this section we construct a trajectory of the separatrix map the projection of which to the domain $\DD$ begins in the $\eps^{1/(8n)}$-neighbourhood of the first endpoint $\varkappa_1$ and finishes in the $\eps^{1/(8n)}$-neighbourhood of the other endpoint $\varkappa_2$. Moreover, the projection lies in the $\eps^{1/(8n)}$-neighbourhood of $\varkappa$.
Afterwards we estimate the time necessary for this passage.

Suppose we have an admissible trajectory $\OO=(\Omega_1,\ldots,\Omega_m)$ with the code $\overline\OO=(\overline\Omega_1,\ldots,\overline\Omega_m)$. By Lemma \ref{lem-main} we construct the next point of the code $\overline \Omega_{m+1}$ by formulas (\ref{newrho})--(\ref{newzeta}). Then by Lemma \ref{lem:r_attach} the extended code corresponds to some extended admissible trajectory of the separatrix map which is close to the code in the metric $d(\cdot,\cdot)$,
(see Definition \ref{dist}). Since
\begin{equation*}
|\rho_m-\overline \rho_m|<\frac{b(2-b)}{b_\rho}, \quad |\zeta_m-\overline\zeta_m|<\frac{b(2-b)\eps^{3/4}}{b_{\zeta}}, \quad |\tau_m-\overline\tau_m|<\frac{b(2-b)}{b_{\tau}},
\end{equation*}
then by (\ref{eq1}) for sufficiently small $\eps$
\begin{equation*}
\quad \big|\mathcal F_{\beta,p}\big(\overline\zeta_{m-1}+\hat\nu(p)\overline t_m\big)-\hat \Theta_{\zeta}(\varepsilon \rho_m,\zeta_m,\tau_m)\big|<\frac{2b(2-b)}{b_\tau}\big\|\Theta\big\|_{C^2}.
\end{equation*}
Therefore, the following dynamical system with discrete time describes dynamics of $\overline \rho_m$ and $\overline \zeta_m$:
\begin{eqnarray}
\label{rho1}
\!\!\!\!\!\!\!\!
\overline \rho_{m+1}-\overline\rho_m \!\! &=& \!\! - v_m+(\overline\tau_{m+1}-\overline\tau_m-t_{m+1})\bH_{\zeta}\big(\eps\overline\rho_m,\overline\zeta_m\big)+r_m,\\
\label{zeta(rho)1}
\!\!\!\!\!\!\!\!
\overline\zeta_{m+1}-\overline\zeta_m\!\! &=& \!\!\hat\nu(\varepsilon \overline \rho_m)\,t_{m+1}+s_m.
\end{eqnarray}
Here
\begin{equation}
\label{v_theta}
v_m = \hat \Omega_{\zeta}(\eps \rho_{m+1},\zeta_m,\tau_m)
\end{equation}
\begin{equation}
\label{t}
\overline t_m - \big (\overline K/\lambda\big)^{1/2}\leq t_m \leq \overline t_m, \qquad t_m \in \mZ,
\end{equation}
\begin{equation}
\big(K_0-\log\varepsilon\big)/\lambda+\big(\overline K/\lambda\big)^{1/2}\leq\overline t_m\leq \big(K_0-\log\varepsilon+\overline K\big)/\lambda,
\end{equation}
\begin{equation}
\label{rs}
|r_m|<\frac{b(2-b)}{b_{\tau}}(1+2\big\|\Theta\big\|_{C^2}), \quad |s_m|=O(\varepsilon^{3/4}|\log \varepsilon|).
\end{equation}
System (\ref{rho1})--(\ref{zeta(rho)1}) may be regarded as a system with a control  $v_m$ and a noise $r_m$ and $s_m$, which satisfy the estimates (\ref{rs}).
According to Lemma  \ref{lem-main}, we choose the control $v_m$ equal to $\mathcal F_{\beta,p}\big(\overline\zeta_{m-1}+\hat\nu(\eps\overline\rho_{m-1})\,\overline t_m\big)$ for an appropriate point $p$; the unit vector $\beta \in \mR^n$ is arbitrary, the integer $\overline t_m$ satisfies inequalities (\ref{overlinetm}).

If we have a trajectory $(\overline \rho_m,\overline \zeta_m)$ of (\ref{rho1})--(\ref{zeta(rho)1}), then according to Lemma \ref{lem-main} there exists a close trajectory $(\rho_m,\zeta_m)$ of the separatrix map. Therefore, it is sufficient construct a trajectory of (\ref{rho1})--(\ref{zeta(rho)1}) with desired properties. More precisely, we prove the following proposition.
\begin{prop}
\label{stm}
Let the function $\Theta$ be typical in the sense of Lemma \ref{lem:typ1}. Then for sufficiently small  $\eps>0$ a trajectory $\{(\overline\rho_i,\zeta_i)\}_{i=1}^M$ of system (\ref{rho1})--(\ref{zeta(rho)1}) exists such that

(i) $|\eps \overline\rho_1-\varkappa_0|<\eps^{1/(8n)}$, $|\eps \overline\rho_M-\varkappa_1|<\eps^{1/(8n)}$,

(ii) $\eps \overline\rho_i \in O(\varkappa, \eps^{1/(8n)})$ for all $i=1,\ldots, M$,

(iii) $M < \frac{C}{\eps |\log\eps|}$.\\
Here $C$ is a constant which does not depend on $\eps$.
\end{prop}
Theorem \ref{main} follows from Proposition \ref{stm}.
Since each step of the separatrix map takes time of the order $O(|\log\eps|)$ for the initial Hamiltonian system, we can estimate the time $T$ of passage through the $\frac{1}{|\log \eps|}$-neighbourhood of the stong resonance:
$$
T=O(M|\log\eps|)=O(\eps^{-1}).
$$

\begin{lem}
\label{lem1}
Let $\pi$ be a segment of type $(a)$ (see Lemma \ref{lem:deform}) with endpoints $\pi_1$ and $\pi_2$. Suppose that $\Omega_0, \ldots \Omega_m$ is an admissible trajectory with the code $\overline\Omega_0, \ldots, \overline\Omega_m$ such that $\eps\overline\rho_m$ lies in the $\eps^{1/(7n)}$-neighbourhood of the point $\pi_1$. Then there exists an admissible trajectory $\Omega'_0, \ldots \Omega'_{m+M}$ with the code $\overline\Omega_0, \ldots, \overline\Omega_{m+M}$ such that

(i) $\eps \overline \rho_{m+1}, \ldots, \eps \overline \rho_{m+M} \in O(\pi,\eps^{1/(8n)})$,

(ii) $|\eps \overline \rho_{m+M} - \pi_2|< \eps^{1/5}$,

(iii) $M\le 3|\pi| (\lambda_{\Theta}\eps)^{-1}$.
\end{lem}

\begin{lem}
\label{lem2}
Let $\pi$ be a segment of type $(b)$ with endpoints $\pi_1$ and $\pi_2$. Suppose that $\Omega_0, \ldots \Omega_m$ is an admissible trajectory with the code $\overline\Omega_0, \ldots, \overline\Omega_m$ such that $\eps\overline\rho_m$ lies in the $\eps^{1/5}$-neighbourhood of the point $\pi_1$. Then there exists an admissible trajectory $\Omega'_0, \ldots \Omega'_{m+M}$ with the code $\overline\Omega_0, \ldots, \overline\Omega_{m+M}$ such that

(i) $\eps \overline \rho_{m+1}, \ldots, \eps \overline \rho_{m+M} \in O(\pi,\eps^{1/(8n)})$,

(ii) $|\eps \overline \rho_{m+M}-\pi_2|<\eps^{1/(7n)}$,

(iii) $M\le 3|\pi|(\lambda_{\Theta} \eps)^{-1}$.
\end{lem}
Proposition \ref{stm} follows from these lemmas. We prove them in the next two subsections.
\begin{rem}
\label{skipbars}
For brevity we skip overlines under $\rho_i$ and $\zeta_i$, always assuming that we deal with a quasi-trajectory.
\end{rem}

\section{Proof of Lemma \ref{lem1}}
\label{Etap1}
Let $\pi$ be a segment of type $(a)$ and
\begin{equation}
\label{string0}
\big|\langle  \nu(\eps \rho_m),h\rangle\big|>\eps^{1/(6n)}, \qquad h \in Z.
\end{equation}
We fix $p=\eps \rho_m$ and introduce the variable $\mu$: $\sqrt\eps\mu=\eps\rho$.
We split $\pi$ into segments of equal lengths $l \in \big[2\eps^{1/5},3\eps^{1/5}\big]$ and denote their endpoints by $x_s$, $1\le s < |\pi| \eps^{-1/5}$.
Suppose we have a trajectory $\{(\sqrt\eps\mu_i,\zeta_i)\}_{i=1}^m=\{(\eps \rho_i,\zeta_i)\}_{i=1}^m$ of system (\ref{rho1})--(\ref{zeta(rho)1}), and $\sqrt\eps\mu_m$ lies in the  $\eps^{1/(7n)}$-neighbourhood of the point $x_1=\pi_1$.
Consider the unit vector
$$
\beta=-\frac{x_2-\sqrt\eps\mu_m}{|x_2-\sqrt\eps\mu_m|}
$$
and the corresponding function $\mathcal F_{\beta,p}$.
Since $\pi$ is a type $(a)$, $\bH_{\zeta}(\eps \rho,\zeta)=0$ in the set (\ref{string0}). Equations  (\ref{rho1})--(\ref{zeta(rho)1}) take the form (see Remark \ref{skipbars}):
\begin{eqnarray}
\label{mu2}
\!\!\!\!\!\!\!\!
\mu_{m+1}-\mu_m \!\! &=& \!\! -\sqrt \varepsilon\big(\mathcal F_{\beta,p}\big(\zeta_{m-1}+\hat\nu(p)\overline t_m\big)+ r_m\big),\\
\label{zeta2}
\!\!\!\!\!\!\!\!
\zeta_{m+1}-\zeta_m\!\! &=& \!\!\hat\nu(p)\,t_{m+1}+s_m,
\end{eqnarray}
here $r_m$, $s_m$, $\overline t_m$ satisfy (\ref{t})--(\ref{rs}).

Let $\overline t$ be an integer number which satisfies the inequalities
\begin{equation}
\label{overlinet}
\big(K_0-\log\varepsilon\big)/\lambda\leq\overline t\leq \big(K_0-\log\varepsilon+\overline K\big)/\lambda-2\big(\overline K/\lambda\big)^{1/2}.
\end{equation}
We put
\begin{equation}
\label{Delta-hatzeta}
\Delta=\hat\nu(p)\Big(\overline t+\Big [\big(\overline K/\lambda \big)^{1/2}\Big]\Big), \quad
\hat \zeta = \zeta_{m-1}+\hat\nu(p)\,\overline t_m.
\end{equation}
\begin{prop}
\label{stm:C}
There exist numbers $\overline t_{m+l}\in \mZ$ such that on the corresponding orbit of (\ref{mu2})--(\ref{zeta2}) the inequalities
\begin{align*}
\overline t \leq \overline t_{m+l} &\leq \overline t+\Big [\big(\overline K/\lambda \big)^{1/2}\Big],\\
\big|\zeta_{m+l-1}+\hat\nu(p)\overline t_{m+l}-\hat \zeta -l \Delta \big|&\leq C l^2 \sqrt \varepsilon \,|\log \varepsilon|
\end{align*}
hold for every $l \in \mZ_+$.
\end{prop}
\textbf{Proof.}
We will prove the proposition by induction in $l$. For $l=0$ the inequalities hold. Suppose they hold for some $l \geq 0$, then from (\ref{mu2}) we obtain
\begin{equation}
\label{mu3}
\mu_{m+l}=\mu_m+O(\sqrt \varepsilon l),
\end{equation}
hence from (\ref{zeta2})
$$
\zeta_{m+l}+\hat\nu(p)\,\overline t_{m+l+1}=\zeta_{m+l-1}+\hat\nu(p)(t_{m+l}+\overline t_{m+l+1})+\psi_l,
$$
where $|\psi_l|<Cl\sqrt\varepsilon|\log\varepsilon|$.
Therefore,
\begin{equation*}
\zeta_{m+l}+\hat\nu(p)\,\overline t_{m+l+1} = \zeta_{m+l-1}+\hat\nu(p)\,\overline t_{m+l}+\hat\nu(p)(t_{m+l}-\overline t_{m+l}+\overline t_{m+l+1})+\psi_l.
\end{equation*}
Since $\overline t_{m+l}-t_{m+l} \in \Big[0,\big(\overline K/\lambda \big)^{1/2}\Big]\cap\mZ$, then we can find
$$
\overline t_{m+l+1} \in \Big[\overline t+\big(\overline K/\lambda \big)^{1/2},\overline t ~+~2\big(\overline K/\lambda \big)^{1/2}\Big]
$$
such that
$$
t_{m+l}-\overline t_{m+l}+\overline t_{m+l+1}=\overline t +\Big [\big(\overline K/\lambda \big)^{1/2}\Big].
$$
Hence
$$
\zeta_{m+l}+\hat\nu(p)\,\overline t_{m+l+1}=\hat \zeta + (l+1) \Delta + \tilde \psi_l, \quad |\tilde\psi_l|\le C(l+l^2)\sqrt\varepsilon|\log\varepsilon|.
$$
\qed

\begin{cor}
From (\ref{mu3}) and (\ref{delta}) we obtain
$$
\zeta_{m+l}+\hat\nu(\sqrt \varepsilon \mu_{m+l})\,\overline t_{m+l+1}=\hat \zeta +(l+1)\Delta+\delta_l, \quad |\delta_l|<2\,Cl^2\sqrt\varepsilon|\log \varepsilon|.
$$
\end{cor}

\begin{theo}{\rm(Dirichlet) \cite{Cas57}}.
Let $\triangle=(\triangle_1,\ldots,\triangle_n)\in \mR^n$, $N \in \mN$. Then there exists $\,q \in \mN$ such that $q<N$ and $\max_i\|q\triangle_i\|_{\mZ}\leq\frac{1}{\sqrt[n]{N}}$. Here $\|\cdot\|_{\mZ}$ is the distance to the closest integer.
\end{theo}

Applying the Dirichlet theorem to $\Delta$ for $N=\eps^{-1/5}$ we obtain a vector $\overline \Delta=\big(\overline\Delta_1,\ldots,\overline\Delta_n\big)$ with rational components with the same denominator $q<N$ and
\begin{equation}
\label{deltas}
\big|\Delta - \overline \Delta\big|_{\infty} \leq 1/(q\sqrt[n]{N}).
\end{equation}
For $q < \eps^{-1/5}$ we put
\begin{align*}
\Sigma= -\sum_{l=0}^{q-1}\mathcal F_{\beta,p}(\zeta_{m+l}+\hat\nu(p)\overline t_{m+l+1}\big),
\end{align*}
\begin{equation}
\Sigma_1 =-\sum_{l=0}^{q-1}\mathcal F_{\beta,p}(\hat \zeta +(l+1)\Delta), \quad \Sigma_2 =-\sum_{l=0}^{q-1} \mathcal F_{\beta,p}\big(\hat \zeta +(l+1)\overline\Delta\,\big).
\end{equation}

Points of the set
\begin{equation*}
\Lambda=\{x=(l+1)\overline \Delta \mod \mZ^n, l=0,\ldots,q-1\}
\end{equation*}
form a lattice on the torus $\mT^n$.
Hence the torus can be split into $q$ equal parallelepipeds $\Pi_l$ whose vertices belong to the lattice $\Lambda$, $\mT^n = \cup_{l-1}^q\Pi_l$, and any set $\Pi_{l'}\cap\Pi_{l''}$, $l' \neq l''$ has empty interior.
Such splitting is not unique. We choose
the fundamental parallelepiped $\Pi$ with the smallest diameter $r_{\min}$.
Thus, we consider the splitting $\{\Pi_l=\Pi+v_l, v_l \in \Lambda\}_{l=0}^{q-1}$ of the torus with $\diam \Pi_l=r_{\min}$.

\begin{dfn}
\label{L-set}
Let $L \subset \{0,\ldots, q-1\}$ be the set of indices for which $\calA(p) \cap \Pi_l = \emptyset$.
\end {dfn}

%$r_{\min}$ is a diameter of the minimal fundamental parallelepiped.

\begin{lem}
\label{res1}
Let the vector $\overline\Delta=\big(\overline \Delta_1,\ldots,\overline \Delta_n\big)$ generate a splitting of the torus $\mT^n$ with a minimal diameter $r_{\min}$. Then there exists a vector $k_*=(\hat k_*,k_{0*}) \in \mZ^{n+1} \setminus \{0\}$ such that $\langle \overline \Delta,\hat k_*\rangle=k_{0*}$ and $|k_*|_{\infty}<c_n/r_{\min}$. Here $c_n$ is a constant, which depends only on $n$.
\end{lem}
Proof of this lemma can be found in Appendix \ref{subsec:dens}.
\begin{lem}
\label{res2}
Let $r>0$ be a constant independent of $\eps$. Then there exists $\eps_0>0$ such that for each $\eps \in (0,\eps_0)$ the integer number $\overline t$ can be found with the following properties:

1) $\overline t$ satisfies inequalities (\ref{overlinet}),

2) $\langle \overline \Delta,\hat k_*\rangle + k_{0*} \neq 0$ for all $k_*=(\hat k_*,k_{0*}) \in \mZ^{n+1}\setminus \{0\}$, $|k_*|_{\infty} < c_n/r$.\\
Here $\overline\Delta$ is a
rational approximation (\ref{deltas}) of $\Delta=\Delta(\overline t)$ that corresponds to the number
$\overline t$ (see (\ref{Delta-hatzeta})).
\end{lem}
\textbf{Proof.}
There exists a number $\eps_0>0$ such that for all $\eps\in(0,\eps_0)$ the set
\begin{equation}
\big\{\mu \colon |\langle \nu(\sqrt\eps\mu),k\rangle| \in (\eps^{1/(6n)},\delta), \dist(\sqrt\eps\mu,\chi)<5\eps^{1/(8n)}\big\}
\end{equation}
does not intersect any resonance strip $S^h_0$ with $|h|_{\infty}<c_n/r$.

Take the smallest number $\overline t_0$ which satisfies the inequality (\ref{overlinet}).
Denote
$$
\Delta_i=\hat\nu(\sqrt \eps \mu_m)\Big(\overline t_0+i+\big[\overline K/\lambda\big]\Big).
$$
Let $\overline \Delta_i$ be a vector with rational components, which corresponds to $\Delta_i$, and let $k_{i} = (\hat k_i,k_{0i}) \in \mZ^{n+1}\setminus\{0\}$ be the shortest integer vector which is orthogonal to $\big(\overline \Delta_i,1\big)$:
$$
\langle \overline\Delta_j,\hat k_j\rangle+k_{0j}=0.
$$
Suppose $\overline t_0+i$, $\overline t_0+j$ also satisfy the inequalities (\ref{overlinet}), $i,j \in \mZ_+$, $i\neq j$. Then
\begin{align*}
\Big|\langle \Delta_i,\hat k_j\rangle+k_{0j}\Big| &\ge \Big | \big |\langle \Delta_i,\hat k_{j}\rangle+k_{0j}-\langle \Delta_j,\hat k_j\rangle-k_{0j} \big|-\big |\langle \Delta_j,\hat k_{j} \rangle + k_{0j} \big| \Big|\\
&\ge\Big|\big\langle(j-i)\hat\nu(\sqrt \eps \mu_m),\hat k_j\big\rangle\Big |- \sqrt n \big|\Delta_j-\overline \Delta_j\big|_{\infty}\big|k_{j}\big |_{\infty}.
\end{align*}
If $|k_j|_{\infty}<c_n/r$ and $\sqrt \eps \mu_m \in O\big(\chi,\eps^{1/(8n)}\big)$, then
$\big|\big\langle\nu(\sqrt \eps \mu_m),\hat k_j\big\rangle\big|>\eps^{1/(6n)}$ for all $\eps\in(0,\eps_0)$.
It follows from the definition of $\eps_0$. Therefore,
$$
\big|\langle \overline \Delta_i,\hat k_j\rangle+k_{0j}\big| \ge \big|\langle \Delta_i,\hat k_j\rangle+k_{0j}\big| - \sqrt n |\Delta_i-\overline\Delta_i|_{\infty}|k_j|_{\infty}>\eps^{1/(6n)}-2\sqrt nc_n\eps^{1/(5n)}/r>0,
$$
and hence the vectors $k_{i}$ and $k_{j}$ are not collinear if $i\ne j$. Since the number of integer vectors $k_{i}$ with $|k_{i}|_{\infty}<c_n/r$ is finite and independent of $\eps$, then there exists a number $i$ such that $\overline t_0+i$ satisfies the inequalities (\ref{overlinet}) and $|k_{i}|_{\infty}>c_n/r$.
\qed
Below we use the number $\overline t=\overline t_0 +i$, for which $|k_i|_{\infty}>c_n/r$, and define $\Delta$ and $\hat \zeta$ by formulas (\ref{Delta-hatzeta}).

\begin{prop}
\label{stm:main}
Consider $\lambda_{\Theta}$ and $\beta$ from Lemma \ref{lem:typ1}.
Then there exists $\eps_0>0$ such that for all $\eps\in(0,\eps_0)$ on the trajectory of the system (\ref{mu2})--(\ref{zeta2})
\begin{equation}
\big|\mu_{m+q+1}-\mu_{m+1}-\sqrt \eps q\lambda_{\Theta}\beta\big|_{\infty}<\sqrt \eps q\lambda_{\Theta}/4.
\end{equation}
\end{prop}
\textbf{Proof.}
Let $F_j=F_j(\zeta)$ be the $j$-th coordinate of $-\mathcal F_{\beta,p}(\zeta)$, $j=1,\ldots,n$. The set of discontinuity of $F_j~\colon~\mT^n~\to~\mR$ is contained in $\calA(p)$ (see (\ref{F}). Let
$$
\overline F_j=\int_{\mT^n}F_j(\zeta)\,d\zeta, \qquad F_{j\max}=\sup_{\mT^n} F_j.
$$
Obviously
$$
0<\overline F_j\leq F_{j\max}\leq\max_{\overline\DD\times\mT^{n+1}}|\Theta_{\zeta}(\eta,\zeta,\tau)|.
$$
The last quantity is finite and independent of $\eps$.
We take $\eps_0>0$ and $r_0>0$ such that
\begin{equation}
\label{A}
\vol\Big(O\big(\calA(p),r_0\big)\Big)<\frac{\overline F_j}{100F_{j\max}}.
\end{equation}
Here $\vol\Big(O\big(\calA(p),r_0\big)\Big)$ is the $n$-dimensional volume of the $r_0$-neighbourhood of the set $\calA(p)$. By the mean-value theorem there exist  $\{\xi_i \in \Pi_i : i \in L\}$ (see Definition \ref{L-set}) such that
$$
\sum_{ i \in L} \int_{\Pi_i}F_j(\zeta)\,d\zeta=\frac{1}{q}\sum_{i \in L}F_j(\xi_i).
$$
Let $C$ be the constant from Proposition  \ref{stm:C}. Note that the functions $F_j$ are Lipschitz continuous on parallelepipeds $\Pi_i$, $i \in L$ with a Lipschitz constant $\|\Theta\|_{C^2}$.
Then, denoting by $\Sigma_j$, $\Sigma_{1j}$ and $\Sigma_{2j}$ the corresponding components of $\Sigma$, $\Sigma_{1}$ and $\Sigma_{2}$, we obtain
\begin{align}
\big|\Sigma_{j}-\Sigma_{1j}\big|&\le \sum_{i \in L}C_jCq^2\sqrt\eps|\log\eps|+2\sum_{i \notin L}F_{j\max}\notag\\
\label{sigma1-sigma2}
&=|L|C_jCq^2\sqrt\eps|\log\eps|+2\,F_{j\max}\big(q-|L|\big),\\
\label{sigma2-sigma3}
\big|\Sigma_{1j}-\Sigma_{2j}\big|&\leq \sum_{i \in L}C_j/\sqrt[n]N+2\sum_{i \notin L}F_{j\max}=|L|C_j/\sqrt[n]N+2\,F_{j\max}\big(q-|L|\big),\\
\label{Teoremaosrednem}
\big|\Sigma_{2j}-q \overline F_j\big|&\leq\sum_{i \in L}C_j \, r_{\min}+2\sum_{i \notin L}F_{j\max}=|L|C_j \,r_{\min}+2\,F_{j\max}\big(q-|L|\big).
\end{align}

We put $r=\min\{r_0,\overline F_j/(100C_j)\}$.
By Lemma \ref{res2} the vector $(\overline \Delta,1)$ is not orthogonal to any non-zero integer vector $k_*$ with $|k_*|_{\infty}<c_n/r$. Then by Lemma \ref{res1} $r_{\min}<\overline F_j/(100C_j)$. Since $1/q<r_{\min}\leq r_0$, we have by (\ref{A})
$$
q-|L|<q\overline F_j/(100F_{j\max}).
$$
Therefore, inequalities (\ref{sigma1-sigma2})--(\ref{Teoremaosrednem}) take the form:
\begin{align}
\label{s1}
&|\Sigma_{j}-\Sigma_{1j}|<q^3C_jC\sqrt\eps|\log\eps|+q\overline F_j/25<q\overline F_j/15,\\
&|\Sigma_{1j}-\Sigma_{2j}|<q/\sqrt[n]N+q\overline F_j/50<q\overline F_j/15,\\
\label{s2}
&|\Sigma_{2j}-q\overline F_j|<qC_j r_{\min}+q\overline F_j/50<q\overline F_j/15.
\end{align}
The first inequality follows from $q<\eps^{-1/5}$.
Then, taking into account (\ref{s1})--(\ref{s2}), we obtain
\begin{align*}
|\Sigma_j-q\overline F_j|&\leq|\Sigma_j-\Sigma_{1j}|+|\Sigma_{1j}-\Sigma_{2j}|+|\Sigma_{2j}-q\overline F_j| < q\overline F_j/5.
\end{align*}
Therefore, $\Sigma_j-q\overline F_j \in \big(-q\overline F_j/5,q\overline F_j/5\big)$. Since
$$
\mu_{m+q+1}-\mu_{m+1}=\sqrt{\varepsilon}\,\Sigma+\sqrt\varepsilon\sum^{q-1}_{l=0}r_{m+l+1}, \quad \Big|\sum^{q-1}_{l=0}r_{m+l+1}\Big|<\frac{2qb(2-b)}{b_{\rho}},
$$
then by Lemma \ref{lem:typ1} we obtain
$$
\big|\mu_{m+q+1}-\mu_{m+1}+\sqrt \eps q\lambda_{\Theta}\beta\big|<\sqrt \eps q\lambda_{\Theta}/4.
$$
\qed
Hence
$$
\big|\sqrt \eps\mu_{m+1}-x_2\big|-\big|\sqrt \eps\mu_{m+q+1}-x_2\big|>3\eps q\lambda_{\Theta}/4,
$$
therefore, the point $\sqrt\eps\mu_{m+q+1}$ is closer than $\sqrt \eps\mu_{m+1}$  to the point $x_2$ at least by a distance $3\eps q\lambda_{\Theta}/4$.

Then we repeat the whole procedure again by defining the map $\mathcal F_{\beta,a}$ for corresponding vector $\beta$ and the point $p=\sqrt \eps \mu_{m+q+1}$. We continue doing it until the extended quasi-trajectory enters the $\eps^{1/5}$-neighbourhood of the point $x_2$. We do the same for the segment $[x_2,x_3]$ etc.

Thus, we obtain the quasi-trajectory that starts in the $\eps^{1/5}$-neighbourhood of $x_1$, enters the $\eps^{1/5}$-neighbourhood of $x_2$ and lies in the $\eps^{1/(8n)}$-neighbourhood of the segment $\pi$.
\begin{cor}
There exists an admissible trajectory of the separatrix map, which enters the $\eps^{1/5}$-neighbourhood of the point $\pi_2$ after $M\leq 3l(\pi)(\lambda_{\Theta}\varepsilon)^{-1}$ steps and lies in the $\eps^{1/(8n)}$-neighbourhood of the segment $\pi$.
\end{cor}
\section{Proof of Lemma \ref{lem2}}
\label{Etap2}
Now we are in the conditions of Lemma \ref{lem2}, $h \in Z$ is the corresponding integer resonant vector, $\pi_1$ and $\pi_2$ are the endpoints of the segment $\pi$, $\eta_0$ is the point of intersection of $\pi$ and $S^h_0$, $\beta$ is a unit vector from Lemma \ref{lem:deform}.

While $\eps \rho_m \in S_{\NN}$ (see (\ref{Ess_and_Nonres})) we extend the quasi-trajectory by the same means as in Section \ref{Etap1}. Therefore, we can consider only the case
\begin{equation}
\label{one-five}
\pi \subset S^h_{\eps^{1/(6n)}}, \quad \dist(\pi_1,\eps \rho_m)<\eps^{1/5}.
\end{equation}
Then to prove Lemma \ref{lem2} it is sufficiently to extend the quasi-trajectory so that it intersects the strip $S^h_{\eps^{1/(6n)}}$ near the segment $\pi$.

Near essential resonances $h$ we have an additional difficulty: the component $(\overline\tau_{m+1}-\tau_m-t_{m+1})\bH_{\zeta}(\eps\rho_m,\zeta_m)$ in system (\ref{rho1}) might become large (see (\ref{H:d})) whenever $\eps \rho_m \in S^h_{\eps^{1/4}}$ and starts affecting the dynamics.
However we will show that in average its contribution is small.
First of all, we obtain more precise expression for the function $\bH_{\zeta}(\cdot,\zeta)$ near $S^h_0$.
The function $\bH_{\zeta}$ in a neighbourhood of a single essential resonance $S_0^{h}$ up to a small error depends only on one resonant phase $\ph = \langle \hat h,\zeta\rangle$. Hence it is natural to introduce the function
\begin{equation*}
\label{H0}
    \hat \bH(\eta,\ph)
  = \sum_{j\in\mZ}
      \phi\Big(j\,
        \frac{\langle h,\nu(\eta)\rangle}
             {\eps^{1/4}}\Big)\,
       H_1^{j h}(\eta)\, e^{2\pi i j\ph}.
\end{equation*}
\begin{lem}
\label{lem:H=H*}
Assume that condition ${\bf H_1 1}$ holds. Then for every $\eta\in U(\pi)$, $\zeta \in \mT^n$
$$
\big|\bH_{\zeta}(\eta,\zeta) - \hat \bH_{\ph}(\eta,\ph)\hat h\big|<C|\log^{-3} \eps|,
$$
$C$ is a constant that does not depend on  $\eps$, $\eta$ or $\zeta$.
\end{lem}
\textbf{Proof.}
Let $j_*=(j,j_0)\in \mZ^{n+1}$.
In the sum
$$
\bH_\zeta(\eta,\zeta)=
\sum_{j_*\in \mZ^{n+1}}2\pi ij\phi \Big( \frac{\langle j_*,\nu(\eta)\rangle } {\eps^{1/4}}\Big) H_1^{j,\,j_0}(\eta)\, e^{2\pi i\langle  j,\zeta\rangle}
$$
all components with $|j_*|_{\infty}<|\log \eps|$ and $j_* \nparallel h$ vanish.
Since $\br > n+3$, the remainder can be estimated as follows:
$$
\Big|\sum_{j_* \nparallel h} 2\pi ij\phi \Big( \frac{ \langle j_*,\nu(\eta)\rangle}{\eps^{1/4}}\Big)
H_1^{j,\,j_0}(\eta)\, e^{2\pi i \langle j,\zeta\rangle}\Big|
\le \sum_{|j_*| \ge |\log \eps|} \frac{ C_H}{|j_*|^{\mathbf r-1}_{\infty}}
\le \frac{C_H c_1}{|\log \eps|^{{\mathbf r}-n-1}} \le C_H c_1|\log^{-3}\eps|,
$$
for some $c_1 = c_1(\br)$. We put $C= C_H c_1$. \qed

Suppose that $\eps \rho_l \in S^h_{\eps^{1/4}}\cap O(\pi,\eps^{1/(8n)})$, $l=m,m+1, \ldots$
We put:
\begin{equation}
\label{mu,1/2}
\hat \nu(\eta_0)=\hat \nu_0, \quad \frac{\partial \hat \nu}{\partial \eta}(\eta_0)=\hat \nu_0', \quad \sqrt\eps\rho_l = \mu_l + \sqrt\eps\,\overline\tau_l\,\mathcal H_{\ph}(\mu_l,\ph_l) \hat h, \quad
\ph_l = \langle \hat h,\zeta_l\rangle,
\end{equation}
$$
\mathcal H(\mu,\ph)
  = \sum_{j\in\mZ}
      \phi\Big(j\,
        \frac{\langle h,\nu(\sqrt \eps \mu)\rangle}
             {\eps^{1/4}}\Big)2\pi ij\,
       H_1^{j h}(\eta_0)\, e^{2\pi i j\ph}.
$$

We split our proof of Lemma  \ref{lem2} in two parts. In the first one we construct a quasi-trajectory that crosses
the strong resonance, i.e. $h = k$. The second part is devoted to other essential resonances in the neighbourhood
of the strong resonance, i.e. $h \nparallel k$. It is easier to consider these cases separately, because we use
different methods to construct the quasi-trajectory for each one. Without loss of generality we can assume that
$\eta_0=0$ in each case.

We use the following lemma in both cases.

\begin{lem}
\label{lem:muzeta3}
For $l = m,m+1,\ldots $ system (\ref{rho1})--(\ref{zeta(rho)1}) can be rewritten as follows:
\begin{eqnarray}
\label{mu-mu}
\!\!\!\!\!\!\!\!
      \mu_{l+1} - \mu_l
 \!\! &=& \!\! -\sqrt\eps v_l
      - \sqrt\eps\,t_{l+1}
         \big(\hat h\mathcal H_{\ph}(\mu_l,\ph_l) - \EE^{(\mu)}_l
         \big), \\
\label{zeta-zeta}
\!\!\!\!\!\!\!\!
      \zeta_{l+1} - \zeta_l
 \!\! &=& \!\!
      \hat \nu_0 t_{l+1}+\sqrt\eps\,t_{l+1}
         \big( \hat \nu'_0 \mu_{l+1} + \EE_l^{(\zeta)} \big) ,
\end{eqnarray}
where
\begin{equation*}
v_l=\mathcal F_{\beta,\eta_*}(\zeta_{l-1}+\hat \nu_*\overline t_l),
\end{equation*}
and the error terms $\EE_l^{(\mu)}$ and $\EE_l^{(\zeta)}$ satisfy the inequalities
\begin{equation}
\label{DIV2}
    t_{l+1} |\EE_l^{(\mu)}|
  < \tilde c |\log \eps|^{-2}, \quad
    |\EE_l^{(\zeta)}| \le \tilde c\,\eps^{1/4}\mu_{l+1}.
\end{equation}
\end{lem}
\textbf{Proof.}
We put in (\ref{rho1})--(\ref{zeta(rho)1}) $m=l$ and make the change of variables (\ref{mu,1/2}). As a result
\begin{eqnarray*}
      \mu_{l+1} - \mu_l
  &=& \sqrt\eps\,v_l
      -\sqrt\eps\,t_{l+1}
        \hat h\mathcal H_{\ph}(\mu_l,\ph_l)
                  + \sum\nolimits_{j=1}^6 E_j, \\
      \zeta_{l+1} - \zeta_l
  &=& \hat \nu_0 t_{l+1}+\big(\hat\nu(\sqrt\eps\mu_{l+1}) - \hat \nu_0 \big)t_{l+1}
                  + \sum\nolimits_{j=7}^{8} E_j.
\end{eqnarray*}
The error terms $E_j$ are estimated as follows:
\begin{eqnarray*}
\!\!\!\!
     E_1
 &=& \sqrt\eps\overline\tau_{l+1}\hat h
       \big(\hat\bH_{\ph}(\eps\rho_l,\ph_l)
             - \mathcal H_{\ph}(\mu_{l+1},\ph_{l+1}) \big)=O(\eps^{3/4} \log^2 \eps), \\
\!\!\!\!
     E_2
 &=& \sqrt\eps\overline\tau_l\hat h
       \big(\mathcal H_{\ph}(\mu_l,\ph_l)-\hat\bH_{\ph}(\eps\rho_l,\ph_l) \big)=O(\eps^{3/2}|\log \eps|), \\
\!\!\!\!
E_3
 &=& \sqrt\eps (\overline \tau_l - \tau_l)\bH_{\zeta}(\eps\rho_l,\zeta_l)=O(\sqrt\eps |\log^{-1} \eps|), \\
\!\!\!\!
     E_4
 &=& \sqrt\eps t_{l+1}\hat h
       \big(\mathcal H_{\ph}(\mu_l,\ph_l) -\hat\bH_{\ph}(\eps\rho_l,\ph_l) \big)=O(\eps^{3/2}\log^2 \eps), \\
 \!\!\!\!
                 E_5
 &=& \sqrt\eps\overline\tau_{l+1}
       \big(\bH_{\zeta}(\eps\rho_l,\zeta_l)-\hat h\hat\bH_{\ph}(\eps\rho_l,\ph_l)\big)+\sqrt\eps\overline\tau_l
       \big(\hat h\hat\bH_{\ph}(\eps\rho_l,\ph_l) - \bH_{\zeta}(\eps\rho_l,\zeta_l)\big) \\
 &+& \sqrt\eps t_{l+1} \big(\hat h\hat\bH_{\ph}(\eps\rho_l,\ph_l)-\bH_{\zeta}(\eps\rho_l,\zeta_l) \big)=O(\sqrt\eps\log^{-2} \eps),\\
\!\!\!\!
E_6
 &=& r_l\sqrt\eps=O\big(\eps^{5/4}|\log\eps|\big), \\
\!\!\!\!
     E_7
 &=& t_{l+1} \big( \hat \nu(\eps\rho_{l}) - \hat \nu(\sqrt\eps\mu_{l+1})
                 \big)=O(\eps|\log\eps|), \\
\!\!\!\!
     E_8
 &=& s_l=O(\eps^{3/4}|\log\eps|).
\end{eqnarray*}
\qed
We assume that the function $\bH$ satisfies the following condition:

\medskip
${\bf H_1 2}$. {\it All critical points $\ph$ of the function $\hat\bH(\eta_*,\cdot)$ are nondegenerate:
\begin{equation}
\label{H12}
\hat \bH_{\ph}(\eta_*,\ph)=0 \quad \text{implies} \quad
\hat \bH_{\ph \ph}(\eta_*,\ph) \ne 0.
\end{equation}}
\begin{rem}
${\bf 1.}$ Recall that $\eta_*$ is an intersection of the polygonal line $\varkappa$ with $S^k_0$.
Then for every $\eta \in S^k_{\delta} \cap O(\chi,\eps^{1/(8n)})$ all critical points of  $\hat\bH(\eta,\cdot)$ are nondegenerate.

${\bf 2.}$ Hypothesis ${\bf H_1 2}$ holds for a generic function $H_1 \in C^{\bf r}(\DD \times \mT^n \times D \times \mT^1)$.
\end{rem}

Our goal is to construct a trajectory of the system (\ref{rho1})--(\ref{zeta(rho)1}) that intersects the strip $S^h_{\eps^{1/(6n)}}$ near the segment $\pi$. We have two cases: $h = k$ and $h \nparallel k$.
\subsection{$h = k$.}
\label{h=k}
In this case $\eta_0=\eta_*$, $\hat \nu_0=\hat\nu(\eta_*)$, $\hat \nu_0'=\hat \nu_*'$.
First we pass through the set
\begin{equation}
\label{Set:1}
\big\{\eta \colon |\langle \nu(\eta),k\rangle|\in(\eps^{1/4},\eps^{1/(6n)}]\big\}.
\end{equation}
It consists of two connected components:
\begin{align}
\label{String1}
R_1&=\big\{\eta \colon \langle \nu(\eta),k\rangle\in(\eps^{1/4},\eps^{1/(6n)}]\big\},\\
R_2&=\big\{\eta \colon \langle \nu(\eta),k\rangle\in[-\eps^{1/(6n)},-\eps^{1/4})\big\}.
\end{align}
We extend the quasi-trajectory so that it intersects $R_1$ (i.e. it approaches $S^k_0$).
Let $\eps  \rho_m$ be the first point of the quasi-trajectory which lies in the strip $S^k_{\eps^{1/(6n)}}\cap R_1$.
The strip $R_2$ is passed similarly.
Then
\begin{equation}
\langle\nu(\eps\rho_m),k\rangle = \eps^{\omega}, \quad \text{for some} \quad \omega \in \big[1/(6n),1/4\big).
\end{equation}
We put
\begin{equation}
\label{mu,omega}
\eps\rho_l = \eps^{\omega}\mu_l, \quad
\ph_l = \langle \zeta_l,\hat k\rangle.
\end{equation}
\begin{prop}
\label{stm:muzeta2}
For every $l = m,m+1,\ldots $ system (\ref{rho1})--(\ref{zeta(rho)1}) can be rewritten as:
\begin{eqnarray}
\label{mu2-mu2}
\!\!\!\!\!\!\!\!
      \mu_{l+1} - \mu_l
 \!\! &=& \!\! -\eps^{1-\omega}\big(v_l
      - t_{l+1}\EE^{(\mu)}_l\big),\\
\label{zeta2-zeta2}
\!\!\!\!\!\!\!\!
      \zeta_{l+1} - \zeta_l
 \!\! &=& \!\!
      \hat \nu_0 t_{l+1}+\eps^{\omega}\,t_{l+1}
         \big(\hat\nu'_* \mu_{l+1} + \EE_l^{(\zeta)} \big) ,
\end{eqnarray}
where
\begin{equation}
\label{vl}
v_l=\mathcal F_{\beta,\eta_*}(\zeta_{l-1}+\hat \nu_*\overline t_l),
\end{equation}
for every unit $\beta$ and $\overline t_l$ such that inequalities (\ref{overlinetm}) hold.
The errors $\EE_l^{(\mu)}$ and $\EE_l^{(\zeta)}$ satisfy
\begin{equation}
\label{DIV}
    |t_{l+1}\EE_l^{(\mu)}| < \tilde c |\log^{-1} \eps|, \quad
    |\EE_l^{(\zeta)}| < \tilde c\,\eps^{\omega}.
\end{equation}
\end{prop}
\textbf{Proof.}
We put in (\ref{rho1})--(\ref{zeta(rho)1}) $m=l$ and make the change of variables (\ref{mu,omega}). As a result
\begin{eqnarray*}
      \mu_{l+1} - \mu_l
  &=& \eps^{1-\omega}\,v_l
                  + E_1, \\
      \zeta_{l+1} - \zeta_l
  &=& \hat \nu_* t_{l+1}+\eps^{\omega}t_{l+1}\hat \nu_*'\mu_{l+1}
                  + E_2+E_3,
\end{eqnarray*}
where the error terms $E_1$, $E_2$ and $E_3$ are estimated as follows:
\begin{eqnarray*}
\!\!\!\!
     |E_1|
 &=& |r_l\eps^{1-\omega}|<\frac{2b(2-b)}{b_\rho}\big\|\Theta\big\|_{C^2}\eps^{1-\omega}, \\
 \!\!\!\!
     E_2
 &=& s_l=O(\eps^{3/4}|\log\eps|), \\
\!\!\!\!
E_3 &=& t_{l+1} \big( \nu(\eps^{\omega}\mu_{l+1}) - \hat \nu_*- \eps^{\omega}\mu_{l+1}\hat \nu_*'
                 \big)=O(\eps^{2\omega}|\log\eps|).
\end{eqnarray*}
\qed
According to item $(3)$ of Lemma \ref{lem:deform} the segment $\pi$ is parallel to a unit vector $\beta \in \mR^n$ such that
\begin{equation}
\label{beta}
\langle \hat \nu'_*\hat k,\beta\rangle>|\hat\nu'_*\hat k|/2
\end{equation}
Consider the following vectors $\mu^\perp_l$, $w_l^\perp$ and the constants $a$, $\alpha_l$ and $w_l$:
\begin{eqnarray}
\label{nota_thetalp2}
& a = \langle \hat \nu'_* \hat k,\beta\rangle, \quad
  \mu_l = \beta \alpha_l + \mu^\perp_l, \qquad
  \mu^\perp_l \perp \hat \nu_*'\hat k, & \\
\nonumber
&   -v_l + t_{l+1} \EE_l^{(\mu)}
  = t_{l+1} (\beta w_l + w_l^\perp), \qquad
    w^\perp_l \perp \hat \nu_*'\hat k.
\end{eqnarray}
Projecting equations (\ref{mu2-mu2}) and (\ref{zeta2-zeta2}) onto the directions $\nu'_* \hat k$ and $\hat k$ respectively, we obtain:
\begin{eqnarray}
\label{alpha2-alpha2}
     \alpha_{l+1} - \alpha_l
 &=& \eps^{1-\omega} \, t_{l+1} w_l, \\
\label{ph2-ph2}
    \ph_{l+1} - \ph_l
 &=& \eps^{\omega}\,t_{l+1}
       \big( a\alpha_{l+1} + \langle \hat k,\EE_l^{(\zeta)}\rangle \big).
\end{eqnarray}
The sequence $\{\alpha_l\}$ describes evolution
of the sequence $\{\eps\overline \rho_l\}$ along $\beta$.
Our next goal is to study the system (\ref{alpha2-alpha2})--(\ref{ph2-ph2}).

\medskip
\textbf{Functions $f$ and $g$.}

Consider piecewise smooth map $G \colon \mT^1 \to \mT^n$, which is
defined in (\ref{G}), with its set of discontinuity  $\Phi$ (see
(\ref{Phi})).

For a small\footnote{but independent of $\eps$} number $\epsilon>0$ let us denote the corresponding neighbourhood of the set $\calA(\eta_*)$ by $O\big(\calA(\eta_*),\epsilon\big)$ (the definition of $\calA(\eta_*)$ you can find in (\ref{F})).

Now define the following functions $f$ and $g$:

$(i)$ For $\zeta \in \mT^n$ such that $\zeta + \hat \nu_*\overline t_{l+1} \in \mT^n\setminus O\big(\calA(\eta_*),\epsilon\big)$, we put
\begin{equation}
\label{f,g}
f_l(\zeta)\eqdef-\frac{\langle\hat\nu'_* \hat k,\mathcal F_{\beta,\eta_*}(\zeta+\hat \nu_*\overline t_{l+1})\rangle}{\overline ta}, \quad g(\ph)\eqdef-\frac{\langle\hat\nu'_* \hat k,G(\ph)\rangle}{\overline ta}.
\end{equation}
The function $g$ is continuous on the set
\begin{equation}
\label{smooth:set}
\{\ph  \in \mT^1 \setminus O\big(\Phi,\epsilon|\hat k|\big)\}
\end{equation}
and has a bounded derivative. According to (\ref{g(ph)-h(ph)})
$$
\sup_{\zeta + \hat \nu_*\overline t_{l+1} \in \mT^n \setminus O(\calA(\eta_*),\epsilon)}\big|g(\langle\zeta,\hat k\rangle)-f_l(\zeta)\big|<\frac{c_{\star}|\hat\nu'_*\hat k|\epsilon_0}{\overline ta}.
$$

$(ii)$
For $\ph \in O(\Phi,\epsilon|\hat k|)$ and $\zeta \in O\big(\calA(\eta_*),\epsilon\big)$ we redefine $f_l$ and $g$ so that the functions $g \colon \mT \to \mR$ and $f_l \colon \mT^n \to \mR$ are continuous and
\begin{equation}
\label{sup}
\sup_{\ph \in \mT}|g'(\ph)|<K(\epsilon,\epsilon_0), \quad \sup_{\zeta \in \mT^n}\big|g(\langle\zeta,\hat k\rangle)-f_l(\zeta)\big|<\frac{c_{\star}|\hat\nu'_*\hat k|\epsilon_0}{\overline ta}\le \frac{2c_{\star}\epsilon_0}{\overline t}.
\end{equation}
Here the constant $K(\epsilon,\epsilon_0)$ depends on $\epsilon$ and $\epsilon_0$ and does not depend on $\eps$.
\begin{lem}
\label{lem:disc}
Let $\kappa$ be a small positive number which does not depend on $\eps$. Then for sufficiently small $\epsilon>0$ there exist integers $\overline t_l$ such that\\
1) $\overline t_l$ satisfy (\ref{overlinetm});\\
2) for every trajectory $(\zeta_l,\mu_l)$, $m\le l \le m+M$ such that
$$
|\langle \zeta_{M+m},\hat k\rangle - \langle \zeta_m,\hat k\rangle| > 1,
$$
the set $O\big(\calA(\eta_*),\epsilon\big)$ contains not more than $\kappa M$ points of the sequence $\{\zeta_l+\hat \nu_*\overline t_l\}$, $m\le l \le m+M$.
\end{lem}
\textbf{Proof.}
The map $f(\zeta)=F_{\beta,\eta_*}(\zeta+\hat \nu_*\overline t_{l+1})$ is discontinuous at the point $\zeta$ in two cases:\\
1) $\langle\zeta+\hat \nu_* \overline t_{l+1},\hat k\rangle \in \hat\Phi$;\\
2) $\tau(\zeta+\hat \nu_*\overline t_{l+1})$ is an endpoint of $\Big[0,\big(\overline K/\lambda\big)^{1/2}\Big]$.

Consider equations (\ref{ph2-ph2}), describing the dynamics of the sequence $\{\ph_l\}=\{\langle \zeta_l,\hat k\rangle\}$.
Its elements monotonically move in the circle with variable step
\begin{equation*}
\Delta_l=\eps^{\omega}\,t_{l+1}\big( a\alpha_{l+1} + \langle \hat
k,\quad \EE_l^{(\zeta)}\rangle \big)<0, \quad d_1\eps^{\omega}|\log\eps|<|\Delta_l|<d_2\eps^{\omega}|\log\eps|,
\end{equation*}
where $d_{1,2}>0$ are independent of $\eps$ and $l$.
For some $M=O(\eps^{-\omega}|\log^{-1}\eps|)$ the sequence  $\{\ph_l\}$, $l=m,m+1,\ldots,m+M$, makes a complete rotation. Therefore, for a sufficiently small $\epsilon>0$ there is a small part of the sequence that lies in the set
$$
O\big(\hat\Phi-\langle\hat \nu_*,\hat k \rangle\overline t_{l+1},\epsilon|\hat k|\big) \subset \mT.
$$
If $\tau(\zeta_l+\hat \nu_*\overline t_{l+1})$ is an endpoint of $\Big[0,\big(\overline K/\lambda\big)^{1/2}\Big]$, then $\tau(\zeta_l) \in \{\overline t_{l+1},\overline t_{l+1}+\big(\overline K/\lambda\big)^{1/2}\}$. Taking $\overline t_{l+1} \pm 1$ instead $\overline t_{l+1}$, we obtain
$\tau(\zeta+\hat \nu_*\overline t_{l+1}) \in \{1,\big(\overline K/\lambda\big)^{1/2}-1\}$.
\qed
\begin{rem}
The numbers $c_{\star}$, $\epsilon$ and $\epsilon_0$ are independent of $\eps$.
\end{rem}
Consider the following discrete system:
\begin{eqnarray}
\label{alpha2.1-alpha2.1}
     \tilde\alpha_{l+1} - \tilde\alpha_l
 &=& \eps^{1-\omega} \, t_{l+1} g(\tilde\ph_l), \\
\label{ph2.1-ph2.1}
    \tilde\ph_{l+1} - \tilde\ph_l
 &=& \eps^{\omega}\,t_{l+1}
       \big( a\tilde\alpha_{l+1} + \langle \hat k,\EE_l^{(\zeta)}\rangle \big).
\end{eqnarray}
%According to the second inequality in (\ref{sup}) and Lemma \ref{lem:disc}
%this system is close to the system (\ref{alpha2-alpha2})--(\ref{ph2-ph2}).

Further arguments are based on the fact that the right-hand side (RHS) of this system is close to the one of system (\ref{alpha2-alpha2})--(\ref{ph2-ph2}), therefore, their solutions with the same initial conditions are close to each other if the (discrete) time variable is not too big.
Below we show that for sufficiently small $\kappa$, $\epsilon$ and $\epsilon_0$ solutions of (\ref{alpha2.1-alpha2.1})--(\ref{ph2.1-ph2.1}) are close to solutions of (\ref{alpha2-alpha2})--(\ref{ph2-ph2}) with the same initial conditions. Thus, it is sufficient to show, that system (\ref{alpha2.1-alpha2.1})--(\ref{ph2.1-ph2.1}) has solutions which cross the resonance.

Divide the equations (\ref{alpha2.1-alpha2.1})--(\ref{ph2.1-ph2.1}) by $\eps^{1-\omega}\,t_{l+1}=O\big(\eps^{1-\omega}|\log\eps|\big)$.
Note, that (\ref{alpha2.1-alpha2.1})--(\ref{ph2.1-ph2.1}) is a
discretization of the continuous system
\begin{equation}
\label{alpha.ph.2.1}
    \dot\alpha
  = g(\ph), \quad
    \dot\ph
  = \eps^{2\omega-1}a\alpha.
\end{equation}
For sufficiently small $\epsilon$ the function $g$ has a negative average:
$$
\overline g=\int_{\mT}g(x)\,dx<0.
$$
The phase space of system (\ref{alpha.ph.2.1}) is the cylinder
\begin{equation}
\label{Zaa2}
\!\!
W=\big\{(\alpha,\ph\mod 1):\alpha \in [-\alpha_{\star},\alpha_{\star}]\big\}.
\end{equation}
Here $\alpha_{\star}>0$ is a constant which does not depend on $\eps$, the circle
$$
W_0=\{(0,\ph)\} \subset W
$$
corresponds to
$\eta_*$.
We associate with  $\overline\Omega_m,\ldots,\overline\Omega_l,\overline\Omega_{l+1},\ldots$ the sequence
\begin{equation}
         (\alpha_m,\ph_m),\ldots,(\alpha_l,\ph_l),
         (\alpha_{l+1},\ph_{l+1}),\ldots
\end{equation}
If $\{\alpha_l\}$ crosses $W_0$, then $\{\eps\overline\rho_l\}$ crosses the resonance $S^k_0$.
Let the point $(\alpha_m,\ph_m) \in W$ correspond to the point
$(\overline\rho_m,\overline\zeta_m,\overline \tau_m,\overline t_m)$ of the quasi-trajectory.
Consider a solution $z(t)=(\alpha(t),\ph(t))$ of (\ref{alpha.ph.2.1}), passing through $(\alpha_m,\ph_m)$. To be more specific, let $z(0)=(\alpha_m,\ph_m)$ and $T=O(\eps^{1-2\omega})$ be the value of time, for which the angle $\ph$ makes exactly one full rotation.

According to the first inequality in (\ref{sup}) the solution $\{\tilde\alpha_l,\tilde\ph_l\}$ of (\ref{alpha2.1-alpha2.1})--(\ref{ph2.1-ph2.1}) with the initial condition  $(\tilde\alpha_m,\tilde\ph_m)=(\alpha_m,\ph_m) \in W$ approximates $z(t)$ with accuracy of the order
$$
O(\eps^{2\omega-1})\times O(\eps^{1-2\omega}|)\times O(\eps^{1-\omega}|\log\eps|)< C_1\eps^{1-\omega}|\log \eps|,
$$
where $C_1>0$ does not depend on $\eps$.
Since the point $z(T)$ is closer than $z(0)$ to the resonance by the quantity of the order $O(\eps^{1-2\omega}|\log \eps^{-1}|)$, therefore, the corresponding to $z(T)$ point $\tilde\alpha_{m+M}$ is also closer to the resonance than $\tilde\alpha_m$ by a quantity of the order
$$
O(\eps^{1-2\omega}|\log^{-1} \eps|)>C_2 \eps^{1-2\omega} |\log^{-1} \eps|> 100C_1\eps^{1-\omega}|\log \eps|.
$$
Here $C_2>0$ is independent of $\eps$.
The number of steps $M$ can be estimated as follows:
$$
M=O(\eps^{1-2\omega})\times O(\eps^{\omega-1}|\log^{-1}\eps|)<C_3\eps^{-\omega}|\log^{-1} \eps|, \quad C_3>0 \mbox{ does not depend on } \eps.
$$
We take
\begin{equation}
\label{kappa,epsilon}
\kappa<\frac{aC_2}{12C_3|\hat\nu'_*\hat k|\sup_{\mT^{n+1}}|\Theta_{\zeta}(\eta_*,\zeta,\tau)|_{\mT^{n+1}}}, \quad \epsilon_0<\frac{aC_2}{12C_3|\hat\nu'_*\hat k|c_{\star}}.
\end{equation}
and apply Lemma \ref{lem:disc}.
Then after $M$ steps the difference between
the solutions $\{(\alpha_l,\ph_l)\}$ and $\{(\tilde\alpha_l,\tilde\ph_l)\}$ of systems (\ref{alpha2-alpha2})--(\ref{ph2-ph2}) and (\ref{alpha2.1-alpha2.1})--(\ref{ph2.1-ph2.1}) respectively with the same initial condition $(\alpha_m,\ph_m)$ is not more than $C_2 \eps^{1-2\omega}|5\log\eps|^{-1}$.
We use inequality $(\ref{DIV})$.
Therefore, the point $\alpha_{m+M}$ is closer than $\alpha_m$ to the resonance $\mT^1$, at least by the quantity $4C_2 \eps^{1-2\omega}|5\log\eps|^{-1}$.
Hence the sequence $\{\alpha_l\}$, $l=m,m+1,\ldots,m+M$ approaches the resonance with an average speed $O(\eps^{1-\omega})$.
Thus the sequence $\{\eps \rho_l\}$, $l=m,m+1,\ldots,m+M$ moves toward the resonance by the quantity of the order $O(\eps^{1-\omega}|\log^{-1} \eps|)$ with average speed $O(\eps)$ per unit of discrete time.

Then we define a new number $\omega \in \big(1/(6n),1/4\big)$ and repeat the above procedure again until we cross $S^k_{\eps^{1/4}}$.
\begin{cor}
The sequence $\{\eps \rho_l\}$, $l=m,\ldots,l+M$, passes through the set (\ref{Set:1}), where $M < D\eps^{1/(6n)-1}$, $D$ is independent of $\eps$. Moreover, the sequence lies in the $\eps^{1/(7n)}$-neighbourhood of the curve $\chi$.
\end{cor}

Now we extend the quasi-trajectory so that it crosses the strip
\begin{equation}
\label{string3}
S^k_{\varepsilon^{1/4}}= \big\{\eta \in \DD \colon \big|\langle  \nu(\eta),k\rangle\big|\le\eps^{1/4}\big\}.
\end{equation}
We use the same method as has been described above. The only difference is that we should take into account the term $(\overline\tau_{l+1}-\tau_l-t_{l+1})\bH_{\zeta}(\eps\rho_l,\zeta_l)$.

Suppose that $(\rho_m,\zeta_m,\overline \tau_m,\overline t_m)$ is the first point of the quasi-trajectory such that $\eps\rho_m \in S^k_{\eps^{1/4}}$. Below we extend the quasi-trajectory, so that $\{\eps\rho_i\}$ crosses $S^k_0$ in the $\eps^{1/(7n)}$-neighbourhood of $\eta_*$.

We define vectors $\mu^\perp_l$, $w_l^\perp$ and scalars $a$, $\alpha_l$ and $w_l$ by (\ref{nota_thetalp2}). Applying Lemma \ref{lem:muzeta3} for $h=k$ and projecting equations (\ref{mu-mu}) and  (\ref{zeta-zeta}) onto the directions  $\nu'_* \hat k$ and $\hat k$ respectively, we obtain
\begin{eqnarray}
\label{alpha3-alpha3}
     \alpha_{l+1} - \alpha_l
 &=& \sqrt\eps \, t_{l+1} \big(w_l - a^{-1}\langle\hat k,\hat\nu'_*\hat k\rangle\mathcal H_{\ph}(\mu_l,\ph_l)\big), \\
\label{ph3-ph3}
    \ph_{l+1} - \ph_l
 &=& \sqrt\eps \,t_{l+1}
       \big( a\alpha_{l+1} + \langle \hat k,\EE_l^{(\zeta)}\rangle \big).
\end{eqnarray}
Consider the following discrete system:
\begin{eqnarray}
\label{alpha3.1-alpha3.1}
     \tilde\alpha_{l+1} - \tilde\alpha_l
 &=& \sqrt\eps \, t_{l+1} \big(g(\tilde\ph_l) - a^{-1}\langle\hat k,\hat \nu'_*\hat k\rangle\mathcal H_{\ph}(\mu_l,\tilde\ph_l)\big), \\
\label{ph3.1-ph3.1}
    \tilde\ph_{l+1} - \tilde\ph_l
 &=& \sqrt\eps \,t_{l+1}
       \big( a\tilde\alpha_{l+1} + \langle \hat k,\EE_l^{(\zeta)}\rangle \big).
\end{eqnarray}
%Here the function $g$ is the same as in the previous subsection.
Divide the equations (\ref{alpha3.1-alpha3.1})--(\ref{ph3.1-ph3.1}) by $\sqrt\eps\,t_{l+1}=O\big(\sqrt\eps|\log\eps|\big)$.
The system (\ref{alpha3.1-alpha3.1})--(\ref{ph3.1-ph3.1}) is a discretization of ODE
\begin{equation}
\label{alpha.ph.}
    \dot\alpha
  =\overline g - \mathcal {\tilde H}_{\ph} (\mu,\ph), \quad
    \dot\ph
  = a\alpha.
\end{equation}
Here
$$
\overline g=\int_{\mT} g(x)\,dx, \quad \mathcal {\tilde H} (\mu,\ph)=a^{-1}\langle\hat k,\hat \nu'_*\hat k\rangle\mathcal H (\mu,\ph)+\int_0^{\ph}(\overline g - g(x))\,dx.
$$
Note, that $\overline g < -\lambda_{\Theta}/(2\overline t)$ and $\mathcal {\tilde H}_{\ph}(\mu,\ph)=\mathcal H_{\ph}(\eta_*,\ph)$ if
\begin{equation}
\label{mu}
|\langle  \nu(\sqrt\eps\mu),k\rangle|<\eps^{1/4}/2.
\end{equation}

Consider dynamics of (\ref{alpha.ph.}) on the cylinder
\begin{equation}
\label{Zaa3}
\!\!
W=\{(\alpha,\ph \!\!\! \mod 1):\alpha \in [\alpha_-,\alpha_+]\},\quad \;
\alpha_-=-c\,\eps^{-1/4},\quad \alpha_+=c\,\eps^{-1/4},
\end{equation}
where $c>0$ does not depend on $\eps$. Denote its boundaries:
\begin{equation*}
\!\!
W_{\pm}=\{(\alpha,\ph) \in W:\alpha=\alpha_{\pm}\}.
\end{equation*}
System (\ref{alpha.ph.}) has the ``energy integral''
$$
E=\frac{a}{2}\alpha^2-\overline g \ph + \mathcal {\tilde H} (\eta_*,\ph).
$$
If $\overline g=0$, then the separatrices of the system (\ref{alpha.ph.}) are doubled. However, if
$\overline g<0$, they split and one of them connects $W_+$ and $W_-$, crossing the resonance
$$
W_0 = \{(0,\ph) \in W\}
$$
at the hyperbolic fixed point $(0,\ph_*)$.
Let
$$
M=\min_{\ph \in \mT}\mathcal {\tilde H} (\eta_*,\ph), \quad A=\max_{\ph \in \mT} |\mathcal {\tilde H}_{\ph \ph} (\eta_*,\ph)|,
$$
\begin{lem}
\label{lem:connect}
Suppose $\ph_0$ satisfies
$$
\mathcal {\tilde H}(\eta_*,\ph_0) \in I(\ph_*)\colon=(M+\overline g^2/(20A), M+\overline g^2/(10A)).
$$
Then there exist $t'$, $t''$ such that
$$
z(t') \in W_+, \quad z(t'') \in W_-, \quad |t'-t''|<2\alpha_+/|\overline g|.
$$
Here $z(t)=(\alpha(t),\ph(t))$ is a trajectory of the system (\ref{alpha.ph.}) such that $z(0)=(0,\ph_0)$.
\end{lem}
Lemma \ref{lem:connect} is proven in \cite{DT}.

Hence there exist trajectories of (\ref{alpha.ph.}) near $\gamma$, connecting the boundaries of $W$.
Assuming that ${\bf H_1 2}$ holds, it is easy to prove the following
\begin{lem}
\label{inthatZ}
Let $\mathcal {\tilde H}(\eta_*,\ph(0)) \in I(\ph_*)$.
Then the first complete rotation of the angle $\ph(t)$ of the corresponding trajectory with the initial condition $(0,\ph(0))$ takes time $$
\tau < c_*|\log \eps|^{1/2}.
$$
Here $c_*$ does not depend on $\eps$.
\end{lem}
Thus, the trajectory $z(t)$, $t \in [0,\tau]$ crosses the cylinder
$$
\hat W = \{(\alpha,\ph) \colon |\alpha|<\overline c |\log \eps|^{-1/2}\}, \quad \overline c^2  = -\overline g|\log \eps|.
$$

Now we are ready to extend the quasi-trajectory in such manner that it crosses $S^k_{\varepsilon^{1/4}}$ near $\eta_*$.
First we pass through the set $W\setminus\hat W$. Both its connected components can be crossed similarly, so we
construct a sequence $\{(\alpha_l,\ph_l)\}$ which approaches the resonance.
Suppose $(\alpha_m,\ph_m) \in W \setminus \hat W$ corresponds to the point
$(\rho_m,\zeta_m,\overline \tau_m,\overline t_m)$. Then
$$
\alpha_m \in [\overline c |\log \eps|^{-1/2}, c\,\eps^{-1/4}].
$$
Consider a solution $z(t)=(\alpha(t),\ph(t))$ of (\ref{alpha.ph.}) passing through $(\alpha_m,\ph_m)$:  $z(0)=(\alpha_m,\ph_m)$, and let $T<2/|a\alpha_m|$ be the time of one complete rotation of the angle $\ph$ of system (\ref{alpha.ph.}). The corresponding discrete solution $\{\tilde\alpha_l,\tilde\ph_l\}$ of (\ref{alpha3.1-alpha3.1})--(\ref{ph3.1-ph3.1}) with the same initial condition $(\tilde\alpha_m,\tilde\ph_m)=(\alpha_m,\ph_m) \in Z$ approximates $z(t)$ with the accuracy
$$
T\times O(\eps|\log\eps|^2)< C_1|\alpha_m|^{-1}\eps|\log\eps|^2,
$$
where $C_1>0$ does not depend on $\eps$.

The point $z(T)$ is closer to the resonance than $z(0)$ by a quantity of the order $O(|\alpha_m\log \eps|^{-1})$. Therefore, the point $\tilde\alpha_{m+M}$, which corresponds to $z(T)$, is closer to the resonance than $\tilde\alpha_m$ by a quantity of the order
$$
O(|\alpha_m\log \eps|^{-1})>C_2|\alpha_m\log \eps|^{-1} > 100C_1|\alpha_m|^{-1}\eps|\log\eps|^2.
$$
Here $C_2>0$ is independent of $\eps$,
$$
M<C_3 \alpha_m^{-1}\eps^{-1/2}|\log \eps|, \quad \mbox{where } \,\,C_3>0 \,\,\mbox{ is independent of} \,\, \eps.
$$
We take $\kappa$ and $\epsilon_0$ which satisfy (\ref{kappa,epsilon}) and apply Lemma \ref{lem:disc}.
Then after $M$ steps the solutions $\{(\alpha_l,\ph_l)\}$ and $\{(\tilde\alpha_l,\tilde\ph_l)\}$ of the discrete systems (\ref{alpha3-alpha3})--(\ref{ph3-ph3}) and (\ref{alpha3.1-alpha3.1})--(\ref{ph3.1-ph3.1}) respectively with the same initial condition differ by a quantity, which does not exceed $1/5C_2|\alpha_m\log \eps|^{-1}$. Therefore, the point $\alpha_{m+M}$ is closer to the resonance than $\alpha_m$ at least by $4/5C_2|\alpha_m\log \eps|^{-1}$.

We continue this procedure for the continuous solution $z(t)$, that passes through $(\alpha_{m+M},\ph_{m+M})$, until the sequence $\{(\alpha_l,\ph_l)\}$ enters the cylinder $\hat W$.
\begin{cor}
There exists an admissible trajectory of the separatrix map the projection of which on $\DD$ crosses $W\setminus\hat W$ with an average speed of the order $O(\eps)$.
\end{cor}

Now we cross the cylinder $\hat W$.
Let $\mathcal {\tilde H}(\eta_*,\ph_m) \in I(\ph_*)$, $|\alpha_m|<\sqrt \eps \log^2 \eps$.
Consider the solution $z(t)$ of (\ref{alpha.ph.}), passing through $(\alpha_m,\ph_m)$.
By Lemma \ref{inthatZ} $z(t)$ crosses $\hat W$ during the time interval $T=O(|\log \eps|^{1/2})$. The corresponding discrete solution $\{(\tilde\alpha_l,\tilde\ph_l)\}$ of (\ref{alpha3.1-alpha3.1})--(\ref{ph3.1-ph3.1}) approximates $z(t)$ with the accuracy
$$
O(|\log \eps|^{1/2})\times O(\eps\log^2\eps)=O(\eps|\log \eps|^{5/2})< C_1\eps|\log \eps|^{5/2}.
$$
Hence, $\{(\tilde\alpha_l,\tilde\ph_l)\}$ crosses $\hat W$ after
$$
M=O(|\log\eps|^{1/2})\times O(\eps^{-1/2}|\log^{-1}\eps|)<C_3\eps^{-1/2}|\log\eps|^{-1/2}
$$
steps. We put
$C_2=2\overline c$,
then for $\kappa$ and $\epsilon_0$ satisfying (\ref{kappa,epsilon}) the inequality
$$
|\alpha_{m+M}-\tilde \alpha_{m+M}| < C_2|\log\eps|^{-1/2}/5
$$
holds. Therefore, the sequence $\{(\alpha_l,\ph_l)\}$ crosses $\hat W$.

\begin{cor}
The corresponding sequence $\{\eps\overline \rho_l\}$, $l=m,\dots,m+M$, crosses (\ref{string3}), for $M<E(k)\eps^{-3/4}$, $E(k)$ is a constant which is independent of $\eps$. Since
(\ref{one-five}) holds and $M \times O(\eps|\log \eps|) < \eps^{1/(7n)}$, then the sequence lies in $\eps^{1/(7n)}$-neighbourhood of $\pi$.
\end{cor}
\begin{cor}
There exists an admissible trajectory of the separatrix map, the projection of which on $\DD$ crosses (\ref{string3}) with an average speed of the order $O(\eps)$ and lies in the $\eps^{1/(8n)}$-neighbourhood of the curve $\chi$.
\end{cor}
\subsection{$h \nparallel k$}
\label{hnparallelk}
In this case
 \begin{equation}
 \label{another res}
 O(\pi,\eps^{1/(7n)})\cap S^{k_*}_{\eps^{1/(7n)}} = \emptyset
 \end{equation}
 for every $k_* \in \mZ^{n+1} \setminus \{0\}$ such that $|k_*|<c_n/r$. Here $c_n$ is the constant from Lemma \ref{res1}, $r$ is from Proposition \ref{stm:main}, therefore, they are independent of $\eps$.
 Thus, we can use the method of extending a quasi-trajectory from Section \ref{Etap1}.
First, we pass through the set
\begin{equation}
\label{Set:2}
\big\{\eta \colon |\langle \nu(\eta),h\rangle|\in(\eps^{1/4},\eps^{1/(6n)}]\big\}.
\end{equation}
It consists of two connected components:
\begin{align}
\label{String2}
R_1&=\big\{\eta \colon \langle \nu(\eta),h\rangle\in(\eps^{1/4},\eps^{1/(6n)}]\big\},\\
R_2&=\big\{\eta \colon \langle \nu(\eta),h\rangle\in[-\eps^{1/(6n)},-\eps^{1/4})\big\}.
\end{align}
Without loss of generality we assume $\eps \rho_m \in R_1$.
We extend the quasi-trajectory so that it intersects $R_1$ (i.e. it  approaches $S^h_0$). The strip $R_2$ is passed similarly.

Suppose we have the quasi-trajectory $\overline\Omega_0,\ldots,\overline\Omega_m$, where
$(\rho_m,\zeta_m,\overline \tau_m,\overline t_m)$ is the first its point such that
$$
\eps \rho_m=\sqrt \eps\mu_m + \eps\,\overline\tau_m\,\mathcal H_{\ph}(\mu_m,\ph_m) \hat h \in R_1.
$$
Note that $\bH_{\zeta}(\eta,\zeta)=0$ if $\eta \in O(\pi,\eps^{1/(7n)})\cap R_1$, and $\eps \overline\rho_m$ is far from strong resonances (condition (\ref{another res})). Therefore, we can use Proposition \ref{stm:main} to extend the quasi-trajectory $\overline \OO$, the projection of which on $\DD$ crosses the strip $R_1$ after $M$ steps, $M\leq 3|\pi|(\lambda_{\Theta}\varepsilon)^{-1}$.

Now we extend the quasi-trajectory which crosses the strip
\begin{equation}
\label{string}
S^h_{\varepsilon^{1/4}}= \big\{\eta \in \DD \colon \big|\langle  \nu(\eta),h\rangle\big|\le\eps^{1/4}\big\}.
\end{equation}
We can apply Lemma \ref{lem:muzeta3}, then the dynamics of the quasi-trajectory is described by (\ref{mu-mu})--(\ref{zeta-zeta}).
By Proposition \ref{stm:main} there exists a natural number $q<\eps^{-1/5}$ such that
\begin{equation}
\label{sum:q}
\Big|-\sum_{i=0}^{q-1}v_i+\sum^{q-1}_{i=0}t_{m+i+1} \EE_{m+i}^{(\mu)}-q\lambda_{\Theta}\beta\Big|<q\lambda_{\Theta}/4,
\end{equation}
where $v_i=\mathcal F_{\beta,\eta_0}\big(\zeta_{m+i}+\hat \nu(\eta_0)\overline t_{m+i+1}\big)$.

Consider the following system:
\begin{eqnarray}
\label{tildemu}
\!\!\!\!\!\!\!\!
      \tilde\mu_{l+1} - \tilde\mu_l
 \!\! &=& \!\! -\sqrt\eps\beta \lambda_{\Theta}\overline t/(2\overline t_{l+1})
      - \sqrt\eps\,t_{l+1}
         \hat h\mathcal H_{\ph}(\tilde\mu_l,\tilde\ph_l), \\
\label{tildezeta}
\!\!\!\!\!\!\!\!
      \tilde\zeta_{l+1} - \tilde\zeta_l
 \!\! &=& \!\!
      \hat \nu_0 t_{l+1}+\sqrt\eps\,t_{l+1} \hat \nu'_0 \tilde\mu_{l+1},
\end{eqnarray}
where $\tilde \ph_l =\langle \tilde\zeta_l, \hat h\rangle$, $l \ge m$ and $\overline t$ satisfies inequalities (\ref{overlinet}).

Let $(\zeta_l,\mu_l)$ and $(\tilde \zeta_l,\tilde \mu_l)$ be trajectories of (\ref{mu-mu})--(\ref{zeta-zeta}) and (\ref{tildemu})--(\ref{tildezeta}) respectively with the same initial conditions $(\zeta_m,\mu_m)=(\tilde \zeta_m,\tilde \mu_m)$.
Then by (\ref{sum:q})
$$
\langle \mu_{m+q}-\mu_{m},\beta \rangle >\langle \tilde \mu_{m+q}-\tilde \mu_{m},\beta \rangle,
$$
i.e. $\sqrt \eps \mu_{m+q}$ moves farther in the direction $\beta$, than  $\sqrt \eps \tilde \mu_{m+q}$.
Thus, it is sufficient to show that trajectories of (\ref{tildemu})--(\ref{tildezeta}) cross the strip $S^h_{\eps^{1/4}}$.
Define a vector $\tilde \mu^\perp_l \in\KK(h)$ and constants $a$, $\tilde \alpha_l$, $w$ by the equations:
\begin{equation}
a = \langle \hat \nu'_0 \hat h,\beta\rangle, \quad \tilde\mu_l = \beta \tilde \alpha_l + \tilde \mu^\perp_l, \quad  \overline  t  w= - \lambda_{\Theta}/2.
\end{equation}
Projecting equations (\ref{tildemu})--(\ref{tildezeta}) on the directions $\hat \nu'_0 \hat h$ and $\hat h$ respectively we obtain:
\begin{eqnarray}
\label{alpha4-alpha4}
     \tilde\alpha_{l+1} - \tilde\alpha_l
 &=& w - a^{-1}\langle\hat h,\hat \nu'_0\hat h\rangle\mathcal H_{\ph}(\tilde\mu_l,\tilde\ph_l)\big), \\
\label{ph4-ph4}
    \tilde\ph_{l+1} - \tilde\ph_l
 &=& \sqrt\eps \,t_{l+1} a\tilde\alpha_{l+1}.
\end{eqnarray}
Equations (\ref{alpha4-alpha4})--(\ref{ph4-ph4}) are discrete analogue of
\begin{equation}
    \dot\alpha
  =w - \mathcal {\tilde H}_{\ph} (\mu,\ph), \quad
    \dot\ph
  = a\alpha.
\end{equation}
We are interested in the dynamics inside the cylinder
\begin{equation}
W=\{(\alpha,\ph \!\!\!\! \mod 1):\alpha \in [\alpha_-,\alpha_+]\},\quad \;
\alpha_-=-c\,\eps^{-1/4},\quad \alpha_+=c\,\eps^{-1/4},
\end{equation}
where $c>0$ is a constant that does not depend on $\eps$.
Here the circle
$$
W_0 = \{(0,\ph) \in W\}.
$$
corresponds to the resonance $S^h_0$.
Applying the same arguments\footnote{This case is even simpler, because we do not have to smoothen the right-hand side.} as in Subsection \ref{h=k} we obtain the existence of the solution of this continuous system with initial conditions on the boundary of cylinder $W$ which crosses $W_0$. The corresponding discrete solution of (\ref{alpha4-alpha4})--(\ref{ph4-ph4}) also crosses the resonance.
\begin{cor}
The sequence $\{\varepsilon\rho_l\}$ $l=m,\ldots,m+M$ intersects the strip $S^h_{\eps^{1/4}}$, where $M=O(\eps^{-3/4})$, and lies in the $\eps^{1/(8n)}$-neighbourhood of the segment $\pi$.
\end{cor}

\section{Appendix}
\subsection{Poincar\'e-Melnikov potential at a strong resonance}
\label{subsec:stres}
We reformulate Lemma \ref{lem:FF} in more invariant terms.
We define
$$
  \mR^n_\alpha = \{p\in\mR^{n+1}: \langle p,\alpha\rangle = 0\}, \quad
  \partial = \grad, \quad
  \partial_v = \langle v,\partial\rangle,
$$
where $\alpha\in\mR^{n+1}$ and $v$ is a vector field.
For any smooth function $\thet : \mT^{n+1}\to\mR$ and a nonzero vector $\alpha\in\mR^{n+1}$ the set
$$
  J = \{x\in\mT^{n+1} : \partial_{\alpha}\thet = 0\}
$$
is a smooth $n$-manifold with singularities only at points, where $\partial\partial_\alpha\thet = 0$.

We are interested in the case
\begin{equation}
\label{cond}
x=(\zeta,\tau), \quad \thet(x)=\Theta(\eta_*,\zeta,\tau), \quad \alpha=(-\nu(\eta_*),1), \quad \kappa \parallel \big(\beta,\langle \nu(\eta_*),\beta\rangle\big).
\end{equation}

Suppose that $\langle k,\alpha\rangle = 0$ for some nonzero
$k\in\mZ^{n+1}$. We define the new angular variable
$\ph = \langle k,x\rangle \bmod 1$, $x\in\mT^{n+1}$. The torus $\mT^{n+1}$ is foliated by the $n$-tori
\begin{equation}
\label{mT}
    \mT^n_\ph
  = \{ x\in\mT^{n+1} :
       \langle k,x\rangle = \ph = \mbox{const}\}.
\end{equation}
Recall (see (\ref{pw_smooth})) that $\WW$ is a set of piecewise smooth maps
\begin{equation*}
  w : \mT^1\to J, \quad
  \ph\mapsto x = w(\ph) \in\mT_\ph^n \cap J.
\end{equation*}
We are interested in the map $\PP \colon \WW\to \mR^n_\alpha$,
\begin{equation}
\label{PP}
  \WW\ni w \mapsto \PP(w) = \int_{\mT^1} \partial\thet(w(\ph))\, d\ph.
\end{equation}

\begin{prop}
\label{prop:FF}
Suppose that $\thet$ is generic. Then there exists
$\varkappa_0 = \varkappa_0(\thet) > 0$ such that for any $\beta\in\mR^n_\alpha$, $|\beta|\in [0,\varkappa_0]$ we have $\PP(w) = \beta$ for some $w\in\WW$.
\end{prop}
The word ``generic'' means that for any
$r = 2,3,\ldots,\infty,\omega$ and $l\in\mZ_+$ there exists an open dense set $U$ in the space $C^r(\mT^{n+1}_x\times B^l_c\to\mR)$ (variables $c$ on the ball $B^l\subset\mR^l$ are regarded as parameters) such that for any $F\in U$ and any $c\in B^l$ the function $\thet = F(\cdot,c)$ satisfies Proposition \ref{prop:FF}.

To prove Lemma \ref{lem:FF} it is sufficiently to apply Proposition \ref{prop:FF} for (\ref{cond}).
Note that the direction of the vector $\beta \in \mR^n$ can be chosen arbitrarily.

\textbf{Proof of Proposition \ref{prop:FF}}.
We construct the map $w$ so that any point $w(\ph)$ is located near a point of maximum or minimum of the function $\thet|_{\mT_\ph^n}$.

${\bf 1}$. Let $x_+(\ph)$ and $x_-(\ph)$ be (global) maximum and minimum of $\thet|_{\mT_\ph^n}$. Then $\partial\thet(x_\pm(\ph))\parallel k$. Therefore,
\begin{equation}
\label{grad=lambda}
    \partial\thet(x_\pm(\ph))
  = \frac{k}{|k|^2}\, \partial_k\thet(x_\pm(\ph)), \qquad
  \ph\in\mT.
\end{equation}
For generic $\thet$ we can assume that the functions
$$
  x_\pm : \mT\to\mT^{n+1} \quad
  \mbox{and} \quad
  \thet\circ x_\pm : \mT\to\mR
$$
are piecewisely smooth and all (except a finite number of) the critical points $x_\pm(\ph)$ are nondegenerate. The functions
$\thet\circ x_\pm : \mT\to\mR$ are continuous.

${\bf 2}$. Let $\ph_1,\ldots,\ph_N\in\mT$ be the set of points at which $x_+$ or $x_-$ is not smooth. We can assume that they are well-ordered on the circle $\mT$. Hence the intervals
$$
  I_1 = (\ph_1,\ph_2), \quad
  I_2 = (\ph_2,\ph_3), \quad
  \ldots, \quad
  I_N = (\ph_N,\ph_1).
$$
form a partition of $\mT$. The functions $x_\pm |_{I_j}$,
$j = 1,\ldots,N$ are smooth and
\begin{equation}
\label{dx/dphi}
  \frac{dx_\pm}{d\ph} - \frac{k}{|k|^2} \perp k, \qquad
  \ph\in\mT^1\setminus\Phi, \quad
  \Phi = \{\ph_1,\ldots,\ph_N\}.
\end{equation}

We consider the set $\WW_0\subset\WW$ such that for any $w\in\WW_0$ there exists a set
\begin{equation}
\label{hat:Phi}
  \hat\Phi = \{\hat\ph_1,\ldots,\hat\ph_N\}, \qquad
  \hat\ph_j\in \bar I_j,
\end{equation}
such that $w$ is smooth on any interval
$$
  I_{j\downarrow} = (\ph_j,\hat\ph_j), \quad
  I_{j\uparrow} = (\hat\ph_j,\ph_{j+1}), \qquad
  j = 1,\ldots,N.
$$
\medskip

Obviously $x_\pm\in\WW_0$ and, moreover, for any set of signs
$\sigma_{j\downarrow},\sigma_{j\uparrow}\in \{+,-\}$ the function
\begin{equation}
\label{XsPh}
  \ph\mapsto x_{\sigma,\Phi}(\ph)
     := x_{\sigma(\ph)}(\ph), \quad
    \sigma(\ph)
  = \left\{ \begin{array}{cc}
       \sigma_{j\downarrow}, &\mbox{ if } \ph\in I_{j\downarrow},\\
       \sigma_{j\uparrow},   &\mbox{ if } \ph\in I_{j\uparrow}.
            \end{array}
    \right.
\end{equation}
also lies in $\WW_0$.

${\bf 3}$. To control the projection of integral (\ref{PP}) to $(\spn k)^\perp$, we use the following

\begin{lem}
$\PP(x_\pm)= 0$.
\end{lem}

{\bf Proof}. By (\ref{grad=lambda}) and (\ref{dx/dphi}) we have:
\begin{eqnarray*}
     \PP(x_\pm)
 &=& \frac{k}{|k|^2} \int_{\mT} \partial_k\thet(x_\pm(\ph))\, d\ph \\
 &=& \int_{\mT} \Big\langle \partial\thet(x_\pm(\ph)),
                  \frac{dx_\pm(\ph)}{d\ph} \Big\rangle \, d\ph
  = \int_{\mT} \frac{d\thet(x_\pm(\ph))}{d\ph}\, d\ph.
\end{eqnarray*}
Since the functions $\thet\circ x_\pm$ are continuous, the last integrals vanish. \qed

\begin{cor}
Image of the map $\PP$ contains zero.
\end{cor}

${\bf 4}$. For generic $\thet$ the functions
$\partial_{k}\thet \circ x_\pm$ are different.

\begin{cor}
For generic $\thet$ there exists $\hat\varkappa = \hat\varkappa(\thet) > 0$ such that the set $\PP(\WW_0)$ contains the interval
$$
  \big\{ \varkappa k / |k| :
         \varkappa \in (-\hat\varkappa,\hat\varkappa) \big\}.
$$
\end{cor}

{\bf 5}. Consider $w\in\WW_0$ in a small neighborhood of $x_{\sigma,\Phi}\in\WW_0$. More precisely, we put
$$
  w = x_{\sigma,\Phi} + \eps v + O(\eps^2),
$$
where $v$ is continuous on $\mT\setminus (\Phi\cup\hat\Phi)$.

For any $\ph\in\mT\setminus (\Phi\cup\hat\Phi)$ condition (\ref{pw_smooth}) implies
\begin{equation}
\label{v}
  \langle k,v\rangle = 0, \quad
  \langle A\alpha,v\rangle = 0,
\end{equation}
where the $(n+1)\times (n+1)$-matrix $A(\ph)$ is the following Hessian:
$$
     A(\ph)
  = \frac{\partial^2\thet}{\partial x^2} \big( x_{\sigma,\Phi}(\ph) \big).
$$
Since $\alpha\perp k$, it follows that equations (\ref{v}) are equivalent to
\begin{equation}
\label{vv}
  \langle k,v\rangle = 0, \quad
  \langle A_k\alpha,v\rangle = 0, \qquad
  A_k = \pi A \pi^*,
\end{equation}
where
$$
  \pi(\cdot) = (\cdot) - \frac{k}{|k|^2} \langle k, \cdot\rangle
$$
is the orthogonal projector to the space $T\mT^n_\ph$ while $\pi^*$ is the restriction to $T\mT^n_\ph$.

Any point $x_{\sigma,\Phi}(\ph)$ is a point of maximum or minimum of $\thet|_{\mT_\ph^n}$. Therefore, any matrix $A_k(\ph)$ is non-negative or non-positive definite. Moreover, for any typical $\thet$ the operators $A_k(\ph)$ are degenerate only for a finite number of values $\ph\in\mT$.
\smallskip

{\bf 6}. We can regard $v$ satisfying (\ref{vv}) as an element of the tangent space $T_{x_{\sigma,\Phi}}\WW$. The component of $\PP$ perpendicular to $k$, equals $\pi\PP$. We have: $\pi\PP(x_{\sigma,\Phi}) = 0$ and
$$
    d\pi\PP(x_{\sigma,\Phi}) v
  = \int_\mT \pi A(\ph) v(\ph)\, d\ph
  = \int_\mT A_k(\ph) v(\ph)\, d\ph.
$$
We see that
$$
    d\pi\PP(x_{\sigma,\Phi}) \big( T_{x_{\sigma,\Phi}}\WW \big)
  = \big\{ p\in\mR^{n+1} :
        \langle p,k\rangle = \langle p,\alpha\rangle = 0 \big\}.
$$
Now Implicit function arguments show that $\PP$-image of a neighborhood of $\WW_0$ in $\WW$ contains a neighborhood of the interval
$\big( -\hat\varkappa k / (2|k|), \hat\varkappa k / (2|k|) \big)$ in
$\mR_\alpha^n$. \qed

\subsection{Multiple resonances}
\label{section:vague}
This subsection contains the proof of the proposition \ref{k,epsilon}.
\begin{lem}
\label{lem:number}
Consider an integer vector $k = (k_1,\ldots,k_{n+1})\in\mZ^{n+1}$ such that\\ $GCD(k_1,\ldots,k_{n+1}) = d_{n+1}$.\footnote{
 By definition for any $l\in\mZ\setminus\{0\}$ we put $GCD(0,l)=|l|$ and $GCD(0,0)=+\infty$}
Then there exist a matrix $\KK \in SL(n+1,\mZ)$ such that the first row coincides with $k$ and
\begin{equation}
\label{KK}
\|\KK\|_2 \le \sqrt{n+1} \, |k|_2.
\end{equation}
\end{lem}
For the proof of the lemma see in \cite{DT} (Lemma 18.1).

Let $x=(x_1,\ldots,x_{n+1})\in \mR^{n+1}$ be a vector with euclidian coordinates, $(\zeta,\tau)=x \mod \mZ^{n+1}$. Recall that
\begin{equation}
\label{mTT}
\mT^n_{\ph}=\{x \in \mR^{n+1} \colon \langle k,x\rangle=\ph \mod \mZ\}.
\end{equation}
Apply a linear transformation to $\mR^{n+1}$ with the matrix $\KK$. Then the images of (\ref{mTT}) will be tori
\begin{equation}
\label{KmT}
\KK\mT^n_{\ph}=\{x \in \mR^{n+1} \colon x_1=\ph\mod \mZ\}.
\end{equation}
By $E=\KK B_a$ we denote the image of the $n$-dimensional ball $B_a \subset \mT^n_{\ph}$. Since (\ref{KK}) and $\det \KK=1$, it follows that the set $E$ contains a $n$-dimensional ball of the radius
$$
a_{\star}=\frac{a}{(n+1)^{n/2}|k|^n_{\infty}}.
$$
Then the proposition \ref{k,epsilon} follows from the following ptoposition (see \cite{DT}):
\begin{prop}
Let the point $\eta\in S^k_0$ be $(a,K)$-vague and
$\epsilon = \epsilon_{a,n,K} < 1 / (2\pi)$, where
\begin{equation*}
    \epsilon_{a,n,K}
  = \frac{4n\log (2l_{a,n} + 1)}{\pi K}, \quad
    l_{a,n}
  = \frac{3\cdot 2^n}{\pi a^n \sin (\pi a/2)} + 1
\end{equation*}
Then there exist $l\in\mZ^{n+1}$, $| l|_\infty \le 2l_{a,n}$ such that $\eta \in S^l_\epsilon$.
\end{prop}
Applying the proposition to the torus
(\ref{KmT}) with $a=a_{\star}$. If the point $\eta\in S^k_0$ is $(a_{\star},K)$-vague, then there exists a resonance  $l \in \mZ^{n+1}$ which is non-parallel to $(1,0,\ldots,0) \in \mZ^{n+1}$.

\subsection{Density of lattices on a torus}
\label{subsec:dens}
Consider the lattice $\mZ^n\subset\mR^n$ in the standard Euclidean space $(\mR^n,\langle\,,\rangle)$. By $\|\cdot\|$ we denote the corresponding norm.
We regard $\mT^n$ as the quotient $\mR^n / \mZ^n$ with the canonical covering map $\pi:\mR^n\to\mT^n$.

We say that the finite set $\Lambda\subset\mT^n$, which contains zero, is a lattice, if for any $x_0,x_1\in\Lambda$ the point $x_2$ satisfying the equation

\begin{equation}
\label{lattice}
  x_2 = \pm x_1 \pm x_0 \bmod\mZ^n,
\end{equation}
also lies in $\Lambda$.

Consider the corresponding lattice in $\mR^n$:
$$
  \hat\Lambda = \pi^{-1}(\Lambda).
$$
We can always find $n$ linearly independent vectors $x_1,\ldots,x_n\in\hat\Lambda$. Such a system of vectors is said to be {\it fundamental} if the parallelepiped $\Pi$ with vertices at the points
\begin{equation}
\label{sigmax}
  \sigma_1 x_1 + \ldots + \sigma_n x_n, \qquad
  \sigma_j\in \{0,1\}.
\end{equation}
contains no other points of $\hat\Lambda$.
Parallelepiped $\Pi$, generated by a fundamental system, is said to be fundamental.

Fundamental system of vectors for a given lattice is obviously not unique.

\begin{lem}
\label{lem:lattice}
Suppose that for any fundamental parallelepiped $\Pi$ its diameter is greater than $\delta > 0$. Then there exists a vector $b\in\mZ^n$,
$0 < |b| < c_* / \delta$
such that
$$
  \langle b,x\rangle = 0 \bmod\mZ \quad
  \mbox{for every $x\in\Lambda$}.
$$
The constant $c_*$ depends only on $n$.
\end{lem}
\textbf{Proof.} First, we show that a fundamental system  $\{x_1,\ldots,x_n\}\subset\mR^n$ can be chosen so that the angles
$\alpha_{i,j} = \mbox{\rm angle}(x_i,x_j)$ are not too acute.
\begin{prop}
\label{prop:Pi}
The fundamental system can be chosen so that $|\cos\alpha_{ij}| \le 1/2$.
\end{prop}
\textbf{Proof}. We choose $x_1,\ldots,x_n$ so that
\begin{eqnarray}
\nonumber
      |x_1|
  &=& \min\{|x| : x\in\Lambda\}, \\
\nonumber
     |x_2|
  &=& \min\{|x| : x\in\Lambda,
             \mbox{ where $x,x_1$ are linearly independent}\},\\
\label{fundamental}
  &\!\ldots\!\!& \\
\nonumber
     |x_n|
  &=& \min\{|x| : x\in\Lambda,
              \mbox{where $x,x_1,\ldots,x_{n-1}$ are linearly independen}\}.
\end{eqnarray}
In particular we have:
\begin{equation}
\label{x<x}
  |x_1| \le |x_2| \le \ldots \le |x_n|.
\end{equation}

\begin{prop}
The system $\{x_1,\ldots,x_n\}$ is fundamental.
\end{prop}

Indeed, suppose that $x\in\Pi\cap\Lambda$ is not a vertex of $\Pi$. Choose the index $j$ such that the system $x_1,\ldots,x_{j-1},x$ is linearly independent while the system $x_1,\ldots,x_j,x$ is linearly dependent. We generate the parallelepiped $\Pi_j$ (of dimension $j$) by the system $x_1,\ldots,x_j$. Let $V_j$ be the set of its vertices. For any vertex $v\in V_j$ we define the $j$-dimensional simplex $S_v\subset\Pi_j$ as the convex hull of $v$ and all vertices from $V_j$ joined with $v$ by an edge. Then
$\Pi_j = \cup_{v\in V_j} S_v$.

Any symplex $S_v$ has an edge issuing from $v$ and having the length $|x_j|$, while by (\ref{x<x}) lengths of other edges issuing from from $v$ are do not exceed $|x_j|$. Therefore, the ball $B_v$ with the center at $v$ an radius $|x_j|$ contains $S_v$. This implies
$\Pi_j \subset \cup_{v\in V_j} B_v$.

We conclude that $x$ belongs to some ball $B_v$. In other words,
$|v-x| < |x_j|$. The point $v-x$ also lies in $\Lambda$ and the vectors
$x_1,\ldots,x_{j-1},v-x$ are linearly independent. We obtain contradiction with definition of $x_j$. This proves the claim.
\medskip

For any $1\le i < j\le n$ we have: $|x_j\pm x_i| \ge |x_j|$. This is equivalent to
$$
  \langle x_j - x_i, x_j - x_i\rangle \ge \langle x_j, x_j\rangle,
$$
which implies
$$
      |\cos\alpha_{ij}|
  \le \frac{|\langle x_i, x_j\rangle|}{|x_i|^2}
  \le \frac12.
$$
\qed

Below we use the fundamental system determined by (\ref{fundamental}).

\begin{cor}
\label{cor:delta<}
By (\ref{x<x}) diameter of $\Pi$ satisfies $\delta < n |x_n|$.
\end{cor}

Translations of $\Pi$ by vectors $x\in\Lambda$ generate covering of $\mT^n$ such that for any two distinct vectors $x_1,x_2\in\Lambda$
$$
  \Int(x_1 + \Pi) \cap \Int(x_2 + \Pi) = \emptyset.
$$
The $(n-1)$-dimensional faces of the parallelepipeds $x+\Pi$ transversal to $x_n$ will be said to be black.

\begin{cor}
\label{cor:h>}
There exists a constant $c$ which depends only on $n$ such that for the constructed above fundamental system $x_1,\ldots,x_n$ the height $h$ of $\Pi$ perpendicular to the black face satisfies
\begin{equation}
\label{h>}
      h \ge c |x_n|.
\end{equation}
\end{cor}

Let $\BB\subset\mT^n$ be the union of all black faces. Then $\BB$ is an $(n-1)$-dimensional (not necessarily smoth) submanifold of $\mT^n$.

Consider the hyperplane
\begin{equation}
\label{TT}
    \hat\TT
  = \spn (x_1,\ldots,x_{n-1})
  = \{\hat x\in\mR^n : \langle p,\hat x\rangle = 0\}.
\end{equation}
We put $\TT = \pi(\hat\TT) \subset \mT^n$. Then $\TT$ is the connected component of $\BB$ that contains $0$.

Let $d_0$ be the maximal diameter of an open ball in $\mT^n\setminus\TT$. The following proposition gives an important information on the vector $p$.

Given a nonzero vector $a\in\mR^n$ let $\hat\TT_a$ be the hyperplane
$$
  \hat\TT_a = \{\hat x\in\mR^n : \langle a,\hat x\rangle = 0\}.
$$
We define $\TT_a = \pi(\hat\TT_a)$.

\begin{prop}
\label{prop:a||b}
Suppose $\mT^n\setminus\TT_a$ contains an open ball $B$ of diameter $d$. Then $a\parallel b\in\mZ^n$, $0 < |b| \le 1/d$.
\end{prop}
{\textbf Proof}. Consider the full preimage
$$
          \pi^{-1}(\TT_a)
      =   \cup_{n\in\mZ^n} (n + \hat\TT_a)
  \subset \mR^n.
$$
Then $\mR^n\setminus\pi^{-1}(\TT_a)$ contains a ball $\hat B$ of diameter $d$.

If $a$ is not parallel to an integer vector, the inequality
$0 < \langle a,m\rangle < \eps$ has a solution $m_\eps\in\mZ^n$ for arbitrarily small $\eps > 0$. The planes $m_\eps k + \hat\TT_a$, $k\in\mZ$ lie in $\pi^{-1}(\TT_a)$ and intersect any ball of diameter
$$
    \frac{\langle m_\eps,a\rangle}{|a|}
  < \frac{\eps}{|a|}.
$$
This means that $a\parallel b$ for some integer $b\ne 0$. We assume that
$GCD(b_1,\ldots,b_n) = 1$. Then
$$
    \min \{\langle m,b\rangle\, :\, m\in\mZ^n,\, \langle m,b\rangle > 0  \}
  = 1.
$$
Let this minimum be taken on $m_*\in\mZ^n$. Then
$$
    \mR^n\setminus\pi^{-1}(\TT_a)
  = \cup_{k\in\mZ} (m_* k + \hat\TT_a).
$$
Hence for any ball $\hat B \subset \mR^n\setminus\pi^{-1}(\TT_a)$ its radius is smaller than $1/b$. \qed

By Proposition \ref{prop:a||b} without loss of generality the vector $p$, see (\ref{TT}), satisfies
\begin{equation}
\label{<d0}
  p\in\mZ^n, \quad
  GCD(p_1,\ldots,p_n) = 1, \quad
  0 < |p| \le 1/d_0.
\end{equation}

Let $s$ be the number of connected components in $\BB$.

\begin{prop}
\label{prop:,<>=0}
(1) $\langle p,x\rangle = 0 \bmod\mZ$ for any $x\in\TT\cap\Lambda$.

(2) $\langle p,x\rangle = 0 \bmod \frac1s\mZ$ for any $x\in\Lambda$.
\end{prop}

{\it Proof of Proposition \ref{prop:,<>=0}}. Assertion (1) follows from (\ref{TT}). Note that
$$
  \BB = \TT\cup (x_n + \TT)\cup\ldots\cup ((s-1) x_n + \TT), \quad
  \TT = (s x_n + \TT).
$$
This implies assertion (2). \qed

\begin{cor}
$\langle sp,x\rangle = 0 \bmod\mZ$ for any $x\in\Lambda$.
\end{cor}

It remains to estimate $|sp|$. We have:
$$
  |sp| < \frac s{d_0} \le \frac{1}{h}
      \le \frac{1}{c_* |x_n|}
      \le \frac{n}{c_*\delta}
$$
where the first inequality follows from (\ref{<d0}), the second one from the definition of $h$ (see Corollary \ref{cor:h>}), the third one from (\ref{h>}) while the forth one from Corollary \ref{cor:delta<}.

\end{document}